\documentclass[11pt]{article}
\usepackage{psfrag,cite,amsmath,graphicx,epsfig,latexsym,subfigure,color}
\usepackage[scriptsize,bf]{caption}
\usepackage{color}
\usepackage{amssymb}
\DeclareMathOperator{\argmin}{\mbox{argmin}}
\def\bn{\hfill \\ \smallskip\noindent}
\def\argmin{\mathop{\rm argmin}}

\def\vx{x}

\def\dom{\mbox{dom\,}}
\def\proj{\mbox{proj}}
\def\prox{\mbox{prox}}
\def\R{\Re}
\def\dist{\mbox{dist\,}}

\newcommand{\beq}{\begin{equation}}
\newcommand{\eeq}{\end{equation}}

\begin{document}

\bigskip
\def\theequation {\thesection.\arabic{equation}}
\def\pn {\par\smallskip\noindent}
\def \bn {\hfill \\ \smallskip\noindent}
\newcommand{\fs}{f_1,\ldots,f_s}
\newcommand{\f}{\vec{f}}
\newcommand{\vecx}{x_1,\ldots,x_m}
\newcommand{\xoy}{x\rightarrow y}
\newcommand{\barx}{{\bar x}}
\newcommand{\bary}{{\bar y}}
\newtheorem{theorem}{Theorem}[section]
\newtheorem{lemma}{Lemma}[section]
\newtheorem{corollary}{Corollary}[section]
\newtheorem{proposition}{Proposition}[section]
\newtheorem{definition}{Definition}[section]
\newtheorem{claim}{Claim}[section]

\newcommand{\newsection}{\setcounter{equation}{0}\section}

\def\br{\break}
\def\smskip{\par\vskip 5 pt}
\def\proof{\bn {\bf Proof.} }
\def\QED{\hfill{\bf Q.E.D.}\smskip}
\def\qed{\quad{\bf q.e.d.}\smskip}

\newcommand{\cM}{\mathcal{M}}

\begin{titlepage}
\title{\bf On the Linear Convergence of the Alternating Direction Method of Multipliers\thanks{The research is supported  by the National Science Foundation, grant number DMS-1015346.}
} \vskip 0.5cm
\author{{ Mingyi Hong\ \ and\ \ Zhi-Quan Luo}\thanks{ Department of Electrical and
Computer Engineering, University of Minnesota, Minneapolis, MN
55455, USA. Email: \texttt{\{mhong, luozq\}@umn.edu} }
\\[20pt]
\emph{Dedicated to the fond memories of a close friend and
collaborator, Paul Y. Tseng}}
\date{August 13, 2012; Revised March 20, 2013}
\maketitle
\vskip 1.0cm

\begin{abstract}

We analyze the convergence rate of the alternating direction method
of multipliers (ADMM) for minimizing the sum of two or more
nonsmooth convex separable functions subject to linear constraints.
Previous analysis of the ADMM typically assumes that the objective
function is the sum of only \emph{two} convex functions defined on
\emph{two} separable blocks of variables even though the algorithm
works well in numerical experiments for three or more blocks.
Moreover, there has been no rate of convergence analysis for the
ADMM without strong convexity in the objective function. In this
paper we establish the global linear convergence of the ADMM for
minimizing the sum of \emph{any} number of convex separable
functions. This result settles a key question regarding the
convergence of the ADMM when the number of blocks is more than two
or if the strong convexity is absent. It also implies the linear
convergence of the ADMM for several contemporary applications
including LASSO, Group LASSO and Sparse Group LASSO without any
strong convexity assumption. Our proof is based on estimating the
distance from a dual feasible solution to the optimal dual solution
set by the norm of a certain proximal residual, and by requiring
 the dual stepsize to be sufficiently small.
\end{abstract}

\vspace*{\fill}

\pn
\noindent {\bf AMS(MOS) Subject Classifications:}  49, 90.

\end{titlepage}


\newsection{\bf Introduction}

Consider the problem of minimizing a separable nonsmooth convex function subject to linear equality constraints:
\begin{equation}\label{eq:1}
\begin{array}{ll}
\mbox{minimize} & f(x)=f_1(x_1)+f_2(x_2)+\cdots+f_K(x_K)\\ [10pt]
\mbox{subject to} & \displaystyle Ex=E_1x_1+E_2x_2+\cdots+E_Kx_K=q\\
& x_k\in X_k,\quad k=1,2,...,K,
\end{array}
\end{equation}
where each $f_k$ is a nonsmooth convex function (possibly with
extended values), $x=(x_1^T,...,x_K^T)^T\in\Re^n$ is a partition of
the optimization variable $x$, $X=\prod_{k=1}^{K}X_k$ is the
feasible set for $x$, and $E=(E_1,E_2,...,E_K)\in\Re^{m\times n}$ is
an appropriate partition of matrix $E$ (consistent with the
partition of $x$) and $q\in\Re^m$ is a vector. Notice that the model
\eqref{eq:1} can easily accommodate general linear inequality
constraints $Ex\ge q$ by adding one extra block. In particular, we
can introduce a slack variable $x_{K+1}\ge0$ and rewrite the
inequality constraint as $Ex-x_{k+1}=q$. The constraint
$x_{K+1}\ge0$ can be enforced by adding a new convex component
function $f_{K+1}(x_{K+1})=i_{\Re^m_+}(x_{K+1})$ to the objective
function $f(x)$, where $i_{\Re^m_+}(x_{K+1})$ is the indicator
function for the nonnegative orthant $\Re^m_+$
\[
i_{\Re^{m}_+}(x_{K+1})=\left\{\begin{array}{ll}
0,& \mbox{if $x_{K+1}\ge 0$ (entry wise),}\\
\infty, & \mbox{otherwise}.
\end{array}
\right.
\]
In this way, the inequality constrained problem with $K$ blocks is reformulated as an equivalent equality constrained convex minimization problem with $K+1$ blocks.

Optimization problems of the form \eqref{eq:1} arise in many emerging applications involving structured convex optimization. For instance, in compressive sensing applications, we are given an observation matrix $A$ and a noisy observation vector $b\approx Ax$. The goal is to estimate the sparse vector $x$ by solving the following $\ell_1$ regularized linear least squares problem:
\[
\begin{array}{ll}
\mbox{minimize}& \|y\|^2+\lambda\|x\|_1\\
\mbox{subject to} &Ax+y=b,
\end{array}
\]
where $\lambda>0$ is a penalty parameter. Clearly, this is a
structured convex optimization problem of the form \eqref{eq:1} with
$K=2$. If the variable $x$ is further constrained to be nonnegative,
then the corresponding compressive sensing problem can be formulated
as a three block ($K=3$) convex separable optimization problem
\eqref{eq:1} by introducing a slack variable. Similarly, in the
stable version of robust principal component analysis (PCA)
\cite{rpca}, we are given an observation matrix $M\in\Re^{m\times
n}$ which is a noise-corrupted sum of a low rank matrix $L$ and a
sparse matrix $S$. The goal is recover $L$ and $S$ by solving the
following nonsmooth convex optimization problem
\[
\begin{array}{ll}
\mbox{minimize}& \|L\|_*+\rho\|S\|_1 +\lambda\|Z\|_F^2\\
\mbox{subject to} &L+S+Z=M
\end{array}
\]
where $\|\cdot\|_*$ denotes the matrix nuclear norm (defined as the sum of the matrix singular eigenvalues), while $\|\cdot\|_1$ and $\|\cdot\|_F$ denote, respectively, the $\ell_1$ and the Frobenius norm of a matrix (equal to the standard $\ell_1$ and $\ell_2$ vector norms when the matrix is viewed as a vector). In the above formulation, $Z$ denotes the noise matrix, and $\rho,\;\lambda$ are some fixed penalty parameters. It is easily seen that the stable robust PCA problem corresponds to the three block case $K=3$ in the problem \eqref{eq:1} with $x=(L,S,Z)$ and
\begin{equation}\label{eq:fs}
f_1(L)=\|L\|_*,\quad f_2(S)=\|S\|_1,\quad f_3(Z)=\|Z\|_F^2,
\end{equation}
while the coupling linear constraint is given $L+S+Z=M$. In image processing applications where the low rank matrix $L$ is additionally constrained to be nonnegative, then the above problem can be reformulated as
\[
\begin{array}{ll}
\mbox{minimize}& \|L\|_*+\rho\|S\|_1 +\lambda\|Z\|_F^2+i_{\Re^{mn}_+}(C)\\
\mbox{subject to} &L+S+Z=M,\ L-C=0,
\end{array}
\]
where $C$ is a slack matrix variable of the same size as $L$, and $i_{\Re^{mn}_+}(\cdot)$ is the indicator function for the nonnegative orthant $\Re_+^{mn}$.
In this case, the stable robust PCA problem is again in the form of \eqref{eq:1}. In particular, it has 4 block variables $(L,S,Z,C)$ and the first three convex functions are the same as in \eqref{eq:fs}, while the fourth convex function is given by $f_4(C)=i_{\Re^{mn}_+}(C)$. The coupling linear constraints are $L+S+Z=M,\ L-C=0$.
Other applications of the form \eqref{eq:1} include the latent variable Gaussian graphical model selection problem, see \cite{latent}.

A popular approach to solving the separable convex optimization problem \eqref{eq:1} is to attach a Lagrange multiplier vector $y$ to the linear constraints $Ex =q$ and add a quadratic penalty, thus obtaining an augmented Lagrangian function of the form
\begin{equation}\label{eq:aug-lagrangian}
L(x;y) =  f(x)+ \langle y, q - Ex \rangle +\frac{\rho}{2}\|q-Ex\|^2,
\end{equation}
where  $\rho\ge0$ is a constant.
The augmented dual function is given by
\begin{equation}\label{1.1.1}
d(y) = \min_{x} \ f(x)+ \langle y, q - Ex \rangle +\frac{\rho}{2}\|q-Ex\|^2
\end{equation}
and the dual problem (equivalent to \eqref{eq:1} under mild conditions) is
\begin{equation}\label{2.3}
\max_y d(y).
\end{equation}
Moreover, if $\rho>0$, then $Ex$ is constant over the set of minimizers of \eqref{1.1.1} (see Lemma~\ref{lm:const-derivative} in Section 2). This implies that the dual function $d(y)$ is differentiable with
\[
\nabla d(y)=q-Ex(y)
\]
where $x(y)$ is a minimizer of \eqref{1.1.1}. Given the differentiability of $d(y)$, it is natural to consider the following dual ascent method to solve the primal problem \eqref{eq:1}
\begin{equation}\label{eq:dual-ascent}
y:=y+\alpha \nabla d(y)=y+\alpha(q-Ex(y)),
\end{equation}
where $\alpha>0$ is a suitably chosen stepsize.
Such a dual ascent strategy is  well suited
for structured convex optimization problems that are amenable to decomposition.  For example, if the objective function $f$ is separable (i.e., of the form given in \eqref{eq:1}) and if we select $\rho=0$, then the minimization in (\ref{1.1.1}) decomposes into $K$ independent minimizations whose solutions frequently can be obtained in a simple form. In addition, the iterations can be implemented in a manner that exploits the sparsity structure of the problem and, in certain network cases, achieve a high degree of parallelism.   Popular choices for the ascent methods include (single) coordinate ascent (see \cite{[3], [6], [7], [19], [14], [21], [28], [29], [32]}), gradient ascent (see \cite{[14], [21], [31]}) and gradient projection \cite{[Gol64],[LeP65]}. (See \cite{[5], [14], [30]} for additional references.)

For large scale optimization problems, it is numerically advantageous to select $\rho>0$. Unfortunately, this also introduces variable coupling in the augmented Lagrangian \eqref{eq:aug-lagrangian}, which makes the exact minimization step in \eqref{1.1.1} no longer decomposable across variable blocks even if $f$ has a separable structure. In this case, it is more economical to minimize \eqref{1.1.1} inexactly by updating the components of $x$ cyclically via the coordinate descent method. In particular, we can apply the Gauss-Seidel strategy to inexactly minimize \eqref{1.1.1}, and then update the multiplier $y$ using an approximate optimal solution of \eqref{1.1.1} in a manner similar to \eqref{eq:dual-ascent}. The resulting algorithm is called the Alternating Direction Method of Multipliers (ADMM) and is summarized as follows (see \cite{Gab76,Ga83,Glow84,Glt89}). In the general context of sums of
monotone operators, the work of \cite{eckstein10} describes a large family of splitting methods for $K\ge 3$ blocks which, when applied to the dual, result in similar but not identical methods to the ADMM algorithm \eqref{eq:admm} below.
\begin{center}
\fbox{
\begin{minipage}{5.2in}
\smallskip
\centerline{\bf Alternating Direction Method of Multipliers (ADMM)}
\smallskip
At each iteration $r\ge 1$, we first update the primal variable blocks in the Gauss-Seidel fashion and then update the dual multiplier using the updated primal variables:
\begin{equation}\label{eq:admm}
\left\{\begin{array}{l}\displaystyle
x_k^{r+1}={\rm arg}\!\min_{x_k\in X_k}L(x_1^{r+1},...,x^{r+1}_{k-1},x_k,x^r_{k+1},...,x^r_K;y^r),\ k=1,2,...,K,\\[10pt]
\displaystyle y^{r+1}=y^r+\alpha(q-Ex^{r+1})=y^r+\alpha\left(q-\sum_{k=1}^KE_kx_k^{r+1}\right),
\end{array}
\right.
\end{equation}
where $\alpha>0$ is the step size for the dual update.
\end{minipage}
}
\end{center}

Notice that if there is only one block ($K=1$), then the ADMM
reduces to the standard augmented Lagrangian method of multipliers
for which the global convergence is well understood (see e.g.,
\cite{bertsekas}). In particular, it is known that, under mild
assumptions on the problem, this type of dual gradient ascent
methods generate a sequence of iterates whose limit points must be
optimal solutions of the original problem (see \cite{[6], [28],
[30]}).  For the special case of ordinary network flow problems, it
is further known that an associated sequence of dual iterates
converges to an optimal solution of the dual (see \cite{[3]}).  The
rate of convergence of dual ascent methods has been studied in the
reference \cite{luo-tseng} which showed that, under mild assumptions
on the problem, the distance to the optimal dual solution set from
any $y \in \Re^m$ near the set is bounded above by the dual optimality
`residual' $\|\nabla d(y)\|$. By using this bound, it can be
shown that a number of ascent methods, including coordinate ascent
methods and a gradient projection method, converge at least linearly
when applied to solve the dual problem (see \cite{[15],[16]}; also
see \cite{[2],[8],[11]} for related analysis).   (Throughout this
paper, by `linear convergence' we mean root--linear convergence
(denoted by R-linear convergence) in the sense of Ortega and
Rheinboldt \cite{[20]}.)

When there are two blocks ($K=2$), the convergence of the ADMM was
studied in the context of  Douglas-Rachford splitting method
\cite{douglas-rachford,eckstein,eckstein-bertsekas} for finding a
zero of the sum of two maximal monotone operators. It is known that
in this case every limit point of the iterates is an optimal
solution of the problem. The recent work of
\cite{goldfarb-ma,goldfarb-ma-scheinberg,HeYuan12} have shown that,
under some additional assumptions, the objective values generated by
the ADMM algorithm and its accelerated version (which performs some
additional line search steps for the dual update) converge at a rate
of $O(1/r)$ and $O(1/r^2)$ respectively. Moreover, if the objective
function $f(x)$ is strongly convex and the constraint matrix $E$ is
row independent, then the ADMM is known to converge linearly to the
unique minimizer of \eqref{eq:1} \cite{LionsMercier}.  [One notable
exception to the strong convexity requirement is in the special case
of linear programming for which the ADMM is linearly convergent
\cite{eckstein}.] More recent convergence rate analysis of the ADMM
still requires at least one of the component functions ($f_1$ or
$f_2$) to be strongly convex and have a Lipschitz continuous
gradient. Under these and additional rank conditions on the
constraint matrix $E$, some linear convergence rate results can be
obtained for a subset of primal and dual variables in the ADMM
algorithm (or its variant); see \cite{Yin,UCLA,Boley}. However, when
there are more than two blocks involved ($K\ge3$), the convergence
(or the rate of convergence) of the ADMM method is unknown, and this
has been a key open question for several decades. The recent work
\cite{mashiqian} describes a list of novel applications of the ADMM
with $K\ge 3$ and motivates strongly for the need to analyze the
convergence of the ADMM in the multi-block case.  The recent
monograph \cite{boyd} contains more details of the history,
convergence analysis and applications of the ADMM and related
methods.

A main contribution of this paper is to establish the global (linear) convergence of the ADMM method for a class of convex objective functions involving any number of blocks ($K$ is arbitrary). The key requirement for the global (linear) convergence is the satisfaction of a certain error bound condition that is similar to that used in the analysis of \cite{luo-tseng}. This error bound estimates the distance from an iterate to the optimal solution set in terms of a certain proximity residual. The class of objective functions that are known to satisfy this error bound condition include many of the compressive sensing applications, such as LASSO \cite{lasso}, Group LASSO \cite{group-lasso} or Sparse Group LASSO \cite{luo-zhang}.

\newsection{ Technical Preliminaries}

Let $f$ be a closed proper convex function in $\Re^n$, let $E$ be an
$m \times n$  matrix,  let $q$ be a vector in $\Re^m$.
%
Let ${\rm dom}\ f$ denote the effective domain of $f$ and let
$\hbox{int}(\hbox{dom } f)$ denote the interior of ${\rm dom}\ f$.
We make the following standing assumptions regarding $f$: \pn {\bf
Assumption A.}
\begin{itemize}
\item [(a)] The global minimum of \eqref{eq:1} is attained and so is its dual optimal value.
The intersection $X\cap \hbox{int}(\hbox{dom } f)\cap \{x\mid
Ex=q\}$ is nonempty.
\item
[(b)] $f=f_1(x_1)+f_2(x_2)+\cdots+f_K(x_K)$, with each $f_k$ further decomposable as
\[
f_k(x_k)=g_k(A_kx_k)+h_k(x_k)
\]
where $g_k$ and $h_k$ are both convex and continuous over their
domains, and $A_k$'s are some given matrices (not necessarily full
column rank, and can be zero).
\item
[(c)] Each $g_k$ is strictly convex and continuously differentiable
on $\hbox{int}(\hbox{dom } g_k)$ with a uniform Lipschitz continuous
gradient
\begin{align}
&\|\nabla A_k^T g_k(Ax_k)-A^T_k\nabla g_k(Ax_k')\|\le
L\|x_k-x_k'\|,~\quad\quad\forall~x_k,x_k'\in X_k\nonumber
\end{align}
where $L>0$ is a constant. 

\item [(d)] Each $h_k$ satisfies either one of the following
conditions
\begin{enumerate}
\item The epigraph of $h_k(x_k)$ is a polyhedral set.
\item $h_k(x_k)=\lambda_k\|x_k\|_1+\sum_{J}w_J\|x_{k,J}\|_2$, where
$x_k=(\cdots, x_{k,J},\cdots)$ is a partition of $x_k$ with $J$
being the partition index.
\item Each $h_k(x_k)$ is the sum of the functions described in the previous
two items.
\end{enumerate}
\item [(e)] For any fixed and finite $y$ and $\xi$, $\sum_k h_k(x_k)$ is finite for all $x\in\{x: L(x;y)\le \xi\}\cap  X$.
\item [(f)] Each submatrix $E_k$ has full column rank.
\item [(g)] The feasible sets $X_k$,
$k=1,\cdots,K$ are compact polyhedral sets.
\end{itemize}

We have the following remarks regarding to the assumptions made.
\begin{enumerate}
\item Each $f_k$ may only contain convex function $h_k$. That is, the strongly convex
part $g_k$ can be absent. Also, since the matrices $A_k$'s are not required to have full column rank,
the overall objective function $f(x)$ is not necessarily strongly convex. In fact, under Assumption A,
the optimization problem \eqref{eq:1} can still have
multiple primal or dual optimal solutions. This makes the convergence (and rate of convergence) analysis of ADMM difficult.
\item Assumption $(e)$ does allow $h_k(\cdot)$ to be an
indicator function, as in this case the set $\{x_k\mid \sum_k
h_k(x_k)=\infty\}\cap X$ is not a subset of $\{x: L(x;y)\le
\xi\}\cap X$ for any given $y$ and $\xi$.

\item Linear term of the form $\langle b_k,x_k\rangle$ is already
included in $h_k$, as its ephigraph is polyhedral. Moreover, from
the assumption that $X_k$ is polyhedral, the feasibility constraint
$x_k\in X_k$ can be absorbed into $h_k$ by adding to it an indicator
function ${i}_{X_k} (x_k)$. To simplify notations,
we will not explicitly write $x_k\in X_k$ in the ADMM update
\eqref{eq:admm} from now on.

\item Assumption (f) is made to ensure that the subproblems for each
$x_k$ is strongly convex. This assumption will be relaxed later when
the subproblems are solved inexactly; see Section \ref{subProximal}.

\item Assumption (g) requires the feasible set of the variables to
be compact, which is needed to ensure that certain error bounds of
the primal and dual problems of \eqref{eq:1} hold. This assumption
is usually satisfied in practical applications (e.g. the consensus
problems) whenever {\it a priori} knowledge on the variable
domain is available. This assumption can be further relaxed; see the
discussion at the end of Section \ref{secMain}.
\end{enumerate}
Under Assumption A, both the primal optimum and the dual optimum values of \eqref{eq:1} are attained and are equal (i.e., the strong duality holds for \eqref{eq:1}) so that
\[
d^*=\max_yL(x;y)=\max_y\left(f(x)+\langle y,q-Ex\rangle + \frac{\rho}{2}\|Ex-q\|^2\right)=\min_{Ex=q}f(x),
\]
where $d^*$ is the optimal value of the dual of \eqref{eq:1}.

Under Assumption A, there may still be multiple optimal solutions for both the primal problem \eqref{eq:1} and its dual problem. We first claim that the dual functional
\begin{equation}\label{eq:dual}
d(y)=\min_xL(x;y)=\min_{x} \ f(x)+ \langle y, q - Ex \rangle +\frac{\rho}{2}\|q-Ex\|^2,
\end{equation}
is differentiable everywhere. Let $X(y)$ denote the set of optimal solutions for \eqref{eq:dual}.
\begin{lemma}\label{lm:const-derivative}
For any $y\in \Re^m$, both $Ex$ and $A_kx_k$, $k=1,2,...,K$, are constant over $X(y)$. Moreover, the dual function $d(y)$ is differentiable everywhere and
\[
\nabla d(y)=q-Ex(y),
\]
where $x(y)\in X(y)$ is any minimizer of \eqref{eq:dual}.
\end{lemma}
\proof
Fix $y\in\Re^m$. We first show that $Ex$ is invariant over $X(y)$. Suppose the contrary, so that there exist two optimal solutions $x$ and $x'$ from $X(y)$ with the property that $Ex\neq Ex'$. Then, we have
\[
d(y)=L(x;y)=L(x';y).
\]
Due to the convexity of $L(x;y)$ with respect to the variable $x$, the solution set $X(y)$ must be convex, implying ${\bar x}=(x+x')/2\in X(y)$. By the convexity of $f(x)$, we have
\[
\frac{1}{2}\left[(f(x)+\langle y,q-Ex\rangle)+(f(x')+\langle y,q-Ex'\rangle)\right]\ge f(\barx)+\langle y,q-E\barx\rangle.
\]
Moreover, by the strict convexity of $\|\cdot\|^2$ and the assumption $Ex\neq Ex'$, we have
\[
\frac{1}{2}\left(\|Ex-q\|^2+\|Ex'-q\|^2\right)>\|E\barx-q\|^2.
\]
Multiplying this inequality by $\rho/2$ and adding it to the previous inequality yields
\[
\frac{1}{2}\left[L(x;y)+L(x';y)\right]>L(\barx;y),
\]
which further implies
\[
d(y)>L(\barx;y).
\]
This contradicts the definition $d(y)=\min_xL(x;y)$. Thus, $Ex$ is
invariant over $X(y)$. Notice that $d(y)$ is a concave function and
its subdifferential is given by \cite{bertsekas}
$$
\partial d(y)=\mbox{Closure of the convex hull}\;\{\;q-Ex(y)\mid x(y)\in X(y)\;\}.
$$
Since $Ex(y)$ is invariant over $X(y)$, the subdifferential $\partial d(y)$ is a singleton. By Danskin's Theorem, this implies that $d(y)$ is differentiable and the gradient is given by $\nabla d(y)=q-Ex(y)$, for any $x(y)\in X(y)$.

A similar argument (and using the strict convexity of $g_k$) shows that $A_kx_k$ is also invariant over $X(y)$. The proof is complete.
\QED

By using Lemma~\ref{lm:const-derivative}, we show below
a Lipschitz continuity property of $\nabla d(y)$, for $y$ over any
level set of $d$. \pn
\begin{lemma}\label{lm:3.2}
Fix any scalar $\eta \le f^*$ and let ${\cal U} = \{\ y\in \Re^m\ |\
d(y) \ge \eta \ \}$.
Then there holds 
\[
 \|\nabla d(y') - \nabla d(y)\|
\le \frac{1}{\rho}\|y' - y\|, \quad \forall\; y' \in {\cal U}, \ y
\in {\cal U}.
\]
\end{lemma}
\proof
Fix any $y$ and $y'$ in ${\cal U}$. Let $x = x(y)$ and $x' = x(y')$
be two minimizers of $L(x;y)$ and $L(x;y')$ respectively. By
convexity, we have
$$z - E^Ty +\rho E^T(Ex-q)= 0\ \ \mbox{ and }\ \
z' - E^Ty'+\rho E^T(Ex'-q)= 0,$$ where $z$ and $z'$ are some
subgradient vectors in the subdifferential $\partial f(x)$ and
$\partial f(x')$ respectively. Thus, we have
\[
\langle z - E^Ty +\rho E^T(Ex-q),x'-x\rangle = 0\] { and }
\[
\langle z' - E^Ty'+\rho E^T(Ex'-q),x-x'\rangle= 0.
\]
Adding the above two equalities yields
\[
\langle z - z' + E^T(y' - y)-\rho E^TE(x'-x) , x' - x \rangle = 0.
\]
Upon rearranging terms and using the convexity property
\[
\langle  z' - z, x' - x \rangle\ge 0,
\]
we get
\[
\langle y' - y , E(x' - x) \rangle = \langle  z' - z, x' - x
\rangle+\rho \|E(x'-x)\|^2 \ge \rho\|E(x' - x)\|^2.
\]
Thus, $\rho \|E(x' - x)\| \le \|y' - y\|$ which together with
$\nabla d(y') - \nabla d(y) = E(x - x')$ (cf.\
Lemma~\ref{lm:const-derivative})  yields
\[
\|\nabla d(y') - \nabla d(y)\|=\|E(x' - x)\| \le {1 \over \rho }\|y
- y'\|.
\]
The proof is complete. \QED

To show the linear convergence of the ADMM method, we need certain
local error bounds around the optimal solution set $X(y)$ as well as
around the dual optimal solution set $Y^*$. To describe these local
error bounds, we first define the notion of a proximity operator.
Let $h:\dom(h)\mapsto \Re$ be a (possibly nonsmooth) convex
function. For every $\vx\in\dom(h)$, the \emph{proximity operator}
of $h$ is defined as \cite{[25]}
\[ 
\prox_{h}(\vx)={\argmin_{u\in\Re^n}}\;\;h(u)+\frac12\|\vx-u\|^2.
\] 
Notice that if $h(\vx)$ is the indicator function of a closed convex set $X$, then
$$
\prox_h(\vx)=\proj_X(\vx),
$$
so the proximity operator is a generalization of the projection operator. In particular, it is known that the proximity operator satisfies the nonexpansiveness property:
\begin{equation}\label{eq:nonexpansiveness}
\|\prox_h(x)-\prox_h(x')\|\le \|x-x'\|,\quad \forall\ x,x'.
\end{equation}

The proximity operator can be used to characterize the optimality condition for a nonsmooth convex optimization problem. Suppose a convex function $f$ is decomposed as $f(x)=g(Ax)+h(x)$ where $g$ is  strongly convex and differentiable, $h$ is a convex (possibly nonsmooth) function, then we can define the \textit{proximal gradient} of $f$ with respect to $h$ as
\[
\tilde \nabla f(\vx):=\vx-\prox_{h}(\vx-\nabla (f(\vx)-h(x)))=\vx-\prox_{h}(\vx-A^T\nabla g(A\vx)).
\]
If $h\equiv 0$, then the proximal gradient $\tilde \nabla f(\vx)=\nabla f(\vx)$.
In general, $\tilde\nabla f(\vx)$ can be used as the (standard) gradient of $f$ for the nonsmooth minimization  $\min_{\vx\in X}f(\vx)$.
For example, $\tilde \nabla f(\vx^*)=0$ iff $\vx^*$ is a global minimizer.

For the Lagrangian minimization problem \eqref{eq:dual} and
under Assumption A,
the work of \cite{tseng09,luo-tseng,luo-zhang} suggests that the size of the proximal gradient
\begin{eqnarray}
\tilde \nabla_xL(x;y)&:=&x-\prox_h\left(x-\nabla_x(L(x;y)-h(x))\right)\nonumber \\
&=&x-\prox_h\left(x-A^T\nabla g(Ax)+E^Ty-\rho E^T(Ex-q)\right)\label{eq:prox-grad}
\end{eqnarray}
can be used to upper bound the distance to the optimal solution set
$X(y)$ of \eqref{eq:dual}. Here $$h(x):=\sum_{k=1}^Kh_k(x_k),\quad
g(Ax):=\sum_{k=1}^Kg_k(A_kx_k)$$ represent the nonsmooth and the
smooth parts of $f(x)$ respectively.

In our analysis of ADMM, we will also need an error bound for the
dual function $d(y)$. Notice that a $y \in\Re^m$ solves (\ref{2.3})
if and only if $ y$ satisfies the system of nonlinear equations
\[
\nabla d(y)=0.
\]
This suggests that the norm of the `residual' $\|\nabla d(y)\|$ may
be a good estimate of  how close $y$ is from solving (\ref{2.3}).
The next lemma says if the nonsmooth part of $f_k$ takes certain
forms, then the distance to the primal and dual optimal solution
sets can indeed be bounded.

\begin{lemma}\label{lm:eb}
Suppose assumptions A(a)---A(e) hold.
\begin{enumerate}
\item If in addition $X$ is a polyhedral set, then
there exists a positive scalar $\tau$ and $\delta$ such that the following error bound holds
\begin{equation}\label{eq:primaleb}
\dist(x,X(y))\le \tau \|\tilde\nabla_x L(x;y)\|,\quad
\end{equation}
for all $(x,y)$ such that $\|\tilde\nabla_x L(x;y)\|\le \delta$, where the proximal gradient $\tilde \nabla_x L(x;y)$ is given by
\eqref{eq:prox-grad}.  Furthermore, if $X$ is also a compact set, then there exists some $\tau>0$ such that
the error bound \eqref{eq:primaleb} holds for all $x\in X\cap {\rm dom}(h)$.
\item
Similarly, if assumption A-(g) also holds, then for any scalar $\zeta$, there exist
positive scalars $\delta$ and $\tau$ such that
\begin{equation}\label{eq:dualeb}
\dist(y,Y^*)=\|y-y^*\|\le \tau\|\nabla d(y)\|,~\mbox{whenever}~
d(y)\ge \zeta~\mbox{and}~\|\nabla d(y)\|\le \delta.
\end{equation}
\end{enumerate}
Moreover, in both cases the constant $\tau$ is \emph{independent} of
the choice of $y$ and $x$.
\end{lemma}




For any fixed $y$, the proof for the first part of Lemma~\ref{lm:eb}
is identical to those of \cite{tseng09,luo-tseng,luo-zhang}, each of
which shows the error bound with different assumptions on the
objective function $f$. In particular, it was shown that
\eqref{eq:primaleb} holds for all $x$ with $\|\tilde\nabla_x
L(x;y)\|\le \delta$ (i.e., sufficiently close to $X(y)$). An
important new ingredient is the claim that the error bound holds
over the compact set $X\cap{\rm dom}(h)$. This can be seen in two
steps as follows: (1) for all $x\in X\cap{\rm dom}(h)$ such that
$\|\tilde\nabla_x L(x;y)\|\le \delta$, the error bound
\eqref{eq:primaleb} is already known to hold; (2) for all $x\in
X\cap {\rm dom}(h)$ such that $\|\tilde\nabla_x L(x;y)\|\ge \delta$,
the ratio
\[
\frac{\dist(x,X(y))}{ \|\tilde\nabla_x L(x;y)\|}
\]
is a continuous function and well defined over the compact set
$X\cap {\rm dom}(h)\cap\left\{x\mid \|\tilde\nabla_x L(x;y)\|\ge
\delta\right\}.$ Thus, the above ratio must be bounded from  above
by a constant $\tau'$ (independent of $y$). Combining (1) and (2)
yields the desired error bound over the set $X\cap {\rm dom}(h)$.

Another new ingredient in Lemma~\ref{lm:eb}
is the additional claim that the constants $\delta$, $\tau$ are both
independent of the choice of $y$. This property follows directly
from a similar property of Hoffman's error bound \cite{Hoffman} (on
which the error bounds of \cite{tseng09,luo-tseng,luo-zhang} are
based) for a feasible linear system $P:=\{x\mid Ax\le b\}$:
\[
\dist(x,P) \le \tau \|[Ax-b]_+\|,\quad \forall\ x\in \Re^n,
\]
where $\tau$ is independent of $b$. In fact, a careful checking of
the proofs of \cite{tseng09,luo-tseng,luo-zhang} shows that the
corresponding error constants $\delta$ and $\tau$ for the augmented
Lagrangian function $L(x;y)$ can be indeed made \emph{independent}
of $y$. We omit the proof of the first part of Lemma~\ref{lm:eb} for
space consideration.

Dual error bounds like the one stated in the second part of the
lemma have been studied previously by Pang \cite{[22]} and by
Mangasarian and Shiau \cite{[17]}, though in different contexts. The
above error bound is `local' in that it holds only for those $y$
that are bounded or near $Y^*$ (i.e., when $\|\nabla d(y)\|\le
\delta$ as opposed to a `global' error bound which would hold for
all $y$ in $\Re^m$). However if in addition $y$ also lies in some
compact set $Y$, then the dual error bound hold true for all $y\in
Y$ (using the same argument as in the preceding paragraph). In the appendix, we include a proof showing that the dual error
bound holds true, for the case where the epigraph of $h_k$ is
polyhedral (which includes $\ell_1$ norm and indicator function for
polyhedral sets). We note that from this proof it is clear that
indeed the value of $\tau$ in the dual error bound does not depend
on the choice of either $x$ or $y$.

Under Assumption A(f), the augmented Lagrangian function $L(x;y)$
(cf.\ \eqref{eq:aug-lagrangian}) is strongly convex with respect to
each subvector $x_k$. As a result, each alternating minimization
iteration of ADMM \eqref{eq:admm}
\[{\color{black}
\vx^{r+1}_k = \argmin_{x_k}
L(\vx^{r+1}_1, ..., \vx^{r+1}_{k-1}, \vx_k, \vx^r_{k+1}, ..., \vx^r_K;y^r), \quad k=1,...,K.
}\]
has a unique optimal solution. Thus the sequence of iterates $\{x^r\}$ of the ADMM are well defined.
%
The following lemma shows that the alternating minimization of the Lagrangian function gives a sufficient descent of the Lagrangian function value.
\begin{lemma}\label{lm:p-descent}
Suppose Assumptions A(b) and A(f) hold.
Then fix any index $r$, we have 
\begin{equation}\label{eq:p-descent}
L(\vx^r;y^r)- L(\vx^{r+1};y^r) \ge \gamma\|\vx^r - \vx^{r+1}\|^2,
\end{equation}
where the constant $\gamma>0$ is independent of $r$ and $y^r$.
\end{lemma}
\proof By assumptions A(b) and A(f) , the augmented Lagrangian
function
$$
L(x;y)=\sum_{k=1}^K\left(f_k(x_k)+\langle y_k,q_k-E_kx_k\rangle\right) +\frac{\rho}{2}\left\|\sum_{k=1}^KE_kx_k-q\right\|^2
$$
is strongly convex in each variable $x_k$ and has a uniform modulus
$\rho \lambda_{\min}(E_k^TE_k)>0$. Here, the notation
$\lambda_{\min}(\cdot)$ denotes the smallest eigenvalue of a
symmetric matrix. This implies that, for each $k$, that
\begin{equation}\label{eq:strongcvx}
L(x;y)-L(x_1,..,x_{k-1},{\bar x}_k,x_{k+1},...,x_K;y)\ge \rho \lambda_{\min}(E_k^TE_k)\|x_k-{\bar x}_k\|^2,
\end{equation}
for all $x$, where $\bar x_k$ is the minimizer of $\min_{x_k}L(x;y)$ (when all other variables $\{x_j\}_{j\neq k}$ are fixed).

Fix any index $r$. For each $k\in\{1,...,K\}$, by ADMM
\eqref{eq:admm}, $x^{r+1}_{k}$ is the minimizer of
$L(x^{r+1}_1,...,x^{r+1}_{k-1},x_k,x^r_{k+1},x^r_{k+2},...,x^r_K;y^r)$.
It follows from  \eqref{eq:strongcvx}
\begin{equation}\label{eq:older}
L(\vx^{r+1}_1, ..., x_{k-1}^{r+1},\vx^r_k, ..., \vx^r_K;y^r)-L(\vx^{r+1}_1, ..., \vx^{r+1}_k, \vx^r_{k+1}, ..., \vx^r_K;y^r)\ge \gamma\|\vx^r_k - \vx^{r+1}_k\|^2,\ \ \forall \; k,
\end{equation}
where
$$\gamma=\rho \min_k\lambda_{\min}(E_k^TE_k)$$
is independent of $r$ and $y^r$.
Summing this over $k$, we obtain the {\color{black} sufficient decrease} condition
\[{\color{black}
L(\vx^r;y^r) - L(\vx^{r+1};y^r) \ge \gamma\|\vx^r - \vx^{r+1}\|^2.
}\]
This completes the proof of Lemma~\ref{lm:p-descent}.

\QED

To prove the linear convergence of the ADMM algorithm, we also need the following lemma which bounds the size of the proximal gradient $\tilde\nabla L(x^r;y^r)$ at an iterate $x^r$.
\begin{lemma}\label{lm:estimate}
Suppose assumptions A(b)---A(c) hold. Let $\{x^r\}$ be generated by
the ADMM algorithm \eqref{eq:admm}. Then there exists some constant
$\sigma>0$ $($independent of $y^r$$)$ such that
\begin{equation}\label{eq:old-estimate}
\|\tilde\nabla L(x^r;y^r)\|\le \sigma\|x^{r+1}-x^r\|
\end{equation}
for all $r\ge1$.
\end{lemma}
\proof
Fix any $r\ge1$ and any $1\le k\le K$. According to the ADMM procedure \eqref{eq:admm}, the variable $x_k$ is updated as follows
\[
x^{r+1}_k=\argmin_{x_k}\left(h_k(x_k)+g_k(A_kx_k)-\langle y^r,E_kx_k\rangle+\frac{\rho}{2}\left\|E_kx_k+\sum_{j<k}E_jx_j^{r+1}+\sum_{j>k}E_jx^r_j-q\right\|^2\right).
\]
The corresponding optimality condition can be written as
\begin{equation}\label{eq:old}
\! x_k^{r+1}=\prox_{h_k}\left[x_k^{r+1}-A_k^T\nabla_{x_k}g_k(A_kx^{r+1}_k)+E_k^Ty^r-\rho E_k^T\left(\sum_{j\le k}E_jx_j^{r+1}+\sum_{j>k}E_jx^r_j-q\right) \right].
\end{equation}
Therefore, we have
\begin{eqnarray}\label{eq:error_norm}
&&\!\!\!\!\!\!\!\!\!\!\!\!\!\!\!\!\!\left\|x_k^{r+1}-\prox_{h_k}\left(x_k^{r}-A_k^T\nabla_{x_k}g_k(A_kx^{r}_k)+E_k^Ty^r-\rho E_k^T\left(Ex^r-q\right)\right)\right\|=\nonumber\\
&&\Bigg\|\prox_{h_k}\left[x_k^{r+1}-A_k^T\nabla_{x_k}g_k(A_kx^{r+1}_k)+E_k^Ty^r+\rho E_k^T\left(\sum_{j\le k}E_jx_j^{r+1}+\sum_{j>k}E_jx^r_j-q\right) \right]\nonumber\\
&&-\prox_{h_k}\left(x_k^{r}-A_k^T\nabla_{x_k}g_k(A_kx^{r}_k)+E_k^Ty^r+\rho E_k^T\left(Ex^r-q\right)\right)\Bigg\|\nonumber\\
&\le&\Big\|(x^{r+1}_k-x_k^r)-A_k^T(\nabla_{x_k}g_k(A_kx^{r+1}_k)-\nabla_{x_k}g_k(A_kx^{r}_k))+\rho E_k^T\sum_{j\le k}E_j(x_j^{r+1}-x_j^r)\Big\|\nonumber\\
&\le &\|x^{r+1}_k-x_k^r\|+L\|A_k^T\|\|A_k\|\|x^{r+1}_k-x_k^r\|+\rho\|E_k^T\|\sum_{j\le k}\|E_j\|\|x_j^{r+1}-x_j^r\|\nonumber\\
&\le &c\|x^{r+1}-x^r\|,\quad \mbox{for some $c>0$ independent of
$y^r$},
\end{eqnarray}
where the first inequality follows from the nonexpansive property of the prox operator \eqref{eq:nonexpansiveness}, and the second inequality is due to the Lipschitz property of the gradient vector $\nabla g_k$ (cf.\ Assumption~A-(c)).
Using this relation and the definition of the proximal gradient $\tilde\nabla L(x^r;y^r)$, we have
\begin{eqnarray*}
\|\tilde\nabla_{x_k} L(x^r;y^r)\|&=& \left\|x^r_k-\prox_{h_k}\left(x_k^{r}-A_k^T\nabla_{x_k}g_k(A_kx_k^r)+E_k^Ty^r-\rho E_k^T\left(Ex^r-q\right)\right)\right\|\\
&\le &\|x^r_k-x^{r+1}_k\|+ \left\|x_k^{r+1}-\prox_{h_k}\left(x_k^{r}-A_k^T\nabla_{x_k}g_k(A_kx_k^r)+E_k^Ty^r-\rho E_k^T\left(Ex^r-q\right)\right)\right\|\\
&\le & (c+1)\|x^{r+1}-x^r\|,\quad \forall \ k=1,2,...,K.
\end{eqnarray*}
This further implies that the entire proximal gradient vector can be bounded by $\|x^{r+1}-x^r\|$:
\[
\|\tilde\nabla L(x^r;y^r)\|\le (c+1)\sqrt{K}\|x^{r+1}-x^r\|.
\]
Setting $\sigma=(c+1)\sqrt{K}$ (which is independent of $y^r$) completes the proof.
\QED

\newsection{Linear Convergence of ADMM}\label{secMain}

Let $d^*$ denote the dual optimal value and $\{x^r,y^r\}$ be the
sequence generated by the ADMM method \eqref{eq:admm}. Due to
assumption A(a), $d^*$ also equals to the primal optimal value.
Further we denote
\begin{equation}\label{eq:dd-gap}
\Delta_d^r=d^*-d(y^r)
\end{equation}
which represents the gap from dual optimality at the $r$-th iteration. The primal gap to optimality at iteration $r$ is defined as
\begin{equation}\label{eq:pp-gap}
\Delta_p^r=L(x^{r+1};y^{r})-d(y^r),\quad r\ge 1.
\end{equation}
Clearly, we have both $\Delta_d^r\ge0$ and $\Delta_p^r\ge0$ for all
$r$. To establish the linear convergence of ADMM, we need several
lemmas to estimate the sizes of the primal and dual optimality gaps
as well as their respective decrease.

Let $X(y^r)$ denote the set of optimal solutions for the following
optimization problem
\[
\min_x L(x;y^r)=\min_xf(x)+\langle y^r, q - Ex \rangle
+\frac{\rho}{2}\|Ex-q\|^2.
\]
We denote
\[
\barx^r=\argmin_{\barx\in X(y^r)}\|\barx-x^r\|.
\]

We first bound the sizes of the dual and primal optimality
gaps.

\begin{lemma}\label{lm:cost-to-go}
Suppose assumptions A(a)---A(e) and A(g) hold. Then for any scalar $\delta>0$, there exists a positive scalar $\tau'$ such that  
\begin{equation}\label{eq:d-gap}
\Delta_d^r\le \tau' \|\nabla d(y^r)\|^2=\tau'\|Ex(y^r)-q\|^2,
\end{equation}
for any $y^r \in \Re^m$ with $\|\nabla d(y^r)\|\le \delta$.
Moreover, there exist positive scalars $\zeta$ and $\zeta^{'}$
$($independent of $y^r)$ such that
\begin{equation}\label{eq:p-gap}
\Delta_p^r\le \zeta\|x^{r+1}-x^r\|^2+\zeta^{'}\|x^r-\bar{x}^r\|^2,\quad\mbox{for all $r\ge1$.}
\end{equation}
\end{lemma}

\proof 
Fix any $y^r$, and let $y^*$ be the optimal dual solution closest to
$y^r$. Then it follows from the mean value theorem that there exists
some $\tilde y$ in the line segment joining $y^r$ and $y^*$ such
that
\begin{eqnarray*}
\Delta_d^r&=&d(y^*)-d(y^r)\\
&=&\langle \nabla d({\tilde y}),y^*-y^r\rangle\\
&=&\langle \nabla d({\tilde y})-\nabla d(y^*),y^*-y^r\rangle\\
&\le &\| \nabla d({\tilde y})-\nabla d(y^*)\|\|y^*-y^r\|\\
&\le & \frac{1}{\rho}\|{\tilde y}-y^*\|\|y^*-y^r\|\\
&\le & \frac{1}{\rho}\|y^r-y^*\|\|y^*-y^r\|\\
&=&\frac{1}{\rho}\|y^*-y^r\|^2
\end{eqnarray*}
where the second inequality follows from Lemma~\ref{lm:3.2}. Recall
from the second part in Lemma \ref{lm:eb} that there exists some
$\tau$ such that
\[
\dist(y^r,Y^*)=\|y^r-y^*\|\le \tau\|\nabla d(y^r)\|.
\]
Combining the above two inequalities yields
\[
\Delta^r_d=d(y^*)-d(y^r)\le \tau'\|\nabla d(y^r)\|^2,
\]
where $\tau'=\tau^2/\rho$ is a constant. This establishes the
bound on the size of dual gap \eqref{eq:d-gap}.

It remains to prove the bound on the primal gap \eqref{eq:p-gap}.
For notational simplicity, let us separate the smooth and nonsmooth
part of the augmented Lagrangian as follows
\begin{equation}\label{eq:Lbar}
L(x;y)=g(x)+h(x)+\langle y,
q-Ex\rangle+\frac{\rho}{2}\|q-Ex\|^2\nonumber:=\bar{L}(x;y)+h(x).
\end{equation}

Let  $x^{r+1}_{k}$ denote the $k$-th subvector of the primal vector
$x^{r+1}$. From the way that the variables are updated \eqref{eq:old}, we have
\begin{align}
x^{r+1}_k&={\rm
prox}_{h_k}\left[x^{r+1}_k-\nabla_{x_k}\bar{L}\left(\{x^{r+1}_{j\le
k}\}, \{x^r_{j}\}_{j>k}; y^r\right)\right]\nonumber\\
&={\rm
prox}_{h_k}\left[x^r_k-\nabla_{x_k}\bar{L}(x^r;y^r)-x^r_k+x^{r+1}_k-\nabla_{x_k}\bar{L}\left(\{x^{r+1}_{j\le
k}\}, \{x^r_{j}\}_{j>k};
y^r\right)+\nabla_{x_k}\bar{L}(x^r;y^r)\right]\nonumber\\
&:= {\rm
prox}_{h_k}\left[x^r_k-\nabla_{x_k}\bar{L}(x^r;y^r)-e^r_k\right]\label{eqAGP}
\end{align}
where the gradient vector $\nabla_{x_k}\bar{L}\left(\{x^{r+1}_{j\le
k}\}, \{x^r_{j}\}_{j>k}; y^r\right)$ can be explicitly expressed as
\[
\nabla_{x_k}\bar{L}\left(\{x^{r+1}_{j\le k}\}, \{x^r_{j}\}_{j>k};
y^r\right)=A^T_k\nabla_{x_k}g(A_k x^{r+1}_k)-E^T_k y^r+\rho
E^T_k\left(\sum_{j\le k}E_j x^{r+1}_j+\sum_{j>k}E_jx^r_j-q\right)\]
and the error vector $e^r_k$ is defined by
\begin{align}
e^r_k:=x^r_k-x^{r+1}_k+\nabla_{x_k}\bar{L}\left(\{x^{r+1}_{j\le
k}\}, \{x^r_{j}\}_{j>k}; y^r\right)-\nabla_{x_k}\bar{L}(x^r;y^r).
\label{eqDefe}
\end{align}

%

Note that we can bound the norm of $e^r_k$ as follows
\begin{align}\label{eq:error_norm_AGP}
\|e^r_k\|&\le
\|x^r_k-x^{r+1}_k\|+\|\nabla_{x_k}\bar{L}\left(\{x^{r+1}_{j\le k}\},
\{x^r_{j}\}_{j>k};
y^r\right)-\nabla_{x_k}\bar{L}(x^r;y^r)\|\nonumber\\
&\le \|x^r_k-x^{r+1}_k\|+\left\|A^T_k\left(\nabla_{x_k}g(A_k
x^{r+1}_k)-\nabla_{x_k}g(A_k x^{r}_k)\right)+\rho
E^T_k\left(\sum_{j\le k}E_j
(x^{r+1}_j-x^r_k)\right)\right\|\nonumber\\
&\le c\|x^r-x^{r+1}\|,
\end{align}
where the constant $c>0$ is independent of $y^r$, and can take the
same value as in \eqref{eq:error_norm}.

Using \eqref{eqAGP}, and by the definition of the proximity
operator, we have the following
\begin{align}
&h_k(x^{r+1}_k)+\langle x^{r+1}_k-x^r_k,
\nabla_{x_k}\bar{L}(x^r;y^r)+e^r_k\rangle+\frac{1}{2}\|x^{r+1}_k-x^r_k\|^2\nonumber\\
&\le h_k(\bar{x}^r_k)+\langle \bar{x}^r_k-x^r_k,
\nabla_{x_k}\bar{L}(x^r;
y^r)+e^r_k\rangle+\frac{1}{2}\|\bar{x}_k^r-x^r_k\|^2.
\end{align}
Summing over all $k=1,\cdots, K$, we obtain
\begin{align*}
&h(x^{r+1})+\langle x^{r+1}-x^r,
\nabla_{x}\bar{L}(x^r;y^r)+e^r\rangle+\frac{1}{2}\|x^{r+1}-x^r\|^2\nonumber\\
&\le h(\bar{x}^r)+\langle \bar{x}^r-x^r, \nabla_{x}\bar{L}(x^r;
y^r)+e^r\rangle+\frac{1}{2}\|\bar{x}^r-x^r\|^2.
\end{align*}
Upon rearranging terms, we obtain
\begin{align}
&h(x^{r+1})-h(\bar{x}^r)+\langle x^{r+1}-\bar{x}^r,
\nabla_{x}\bar{L}(x^r;y^r)\rangle\le
\frac{1}{2}\|\bar{x}^r-x^r\|^2-\langle x^{r+1}-\bar{x}^r,
e^r\rangle\label{eqBoundDifference}.
\end{align}

Also, we have from the mean value theorem that there exists some
$\tilde{x}$ in the line segment joining $x^{r+1}$ and $\bar{x}^r$
such that
\[\bar{L}(x^{r+1};
y^r)-\bar{L}(\bar{x}^r; y^r)=\langle \nabla_x \bar{L}(\tilde{x};
y^r), x^{r+1}-\bar{x}^r\rangle. \]

Using the above results, we can bound $\Delta^r_p$ by
\begin{align}
\Delta^r_p&=L(x^{r+1}; y^r)-L(\bar{x}^r; y^r)\nonumber\\
&=\bar{L}(x^{r+1}; y^r)-\bar{L}(\bar{x}^r;
y^r)+h(x^{r+1})-h(\bar{x}^r)\nonumber\\
&=\langle \nabla_x \bar{L}(\tilde{x}; y^r),
x^{r+1}-\bar{x}^r\rangle+h(x^{r+1})-h(\bar{x}^r)\nonumber\\
&=\langle \nabla_x \bar{L}(\tilde{x}; y^r)-\nabla_x \bar{L}({x}^r;
y^r), x^{r+1}-\bar{x}\rangle+\langle \nabla_x \bar{L}({x}^r; y^r),
x^{r+1}-\bar{x}\rangle+h(x^{r+1})-h(\bar{x}^r)\nonumber\\
&\le\langle \nabla_x \bar{L}(\tilde{x}; y^r)-\nabla_x \bar{L}({x}^r;
y^r), x^{r+1}-\bar{x}\rangle+\frac{1}{2}\|\bar{x}^r-x^r\|^2+c\sqrt{K}\|x^{r+1}-x^r\|\|x^{r+1}-\bar{x}^r\| \quad \quad  \nonumber\\
&\le
\left(\sum_{k=1}^{K}L\|A_k\|^T\|A_k\|+\rho\|E^T E\|\right)\|\tilde{x}-x^r\|\|x^{r+1}-\bar{x}^r\|\nonumber\\
&\quad\quad\quad+\frac{1}{2}\|\bar{x}^r-x^r\|^2+c\sqrt{K}\|x^{r+1}-x^r\|\|x^{r+1}-\bar{x}^r\|\nonumber\\
&\le \left(\sum_{k=1}^{K}L\|A_k\|^T\|A_k\|+\rho\|E^T
E\|\right)\left(\|x^{r+1}-x^r\|+\|\bar{x}^{r}-x^r\|\right)^2\nonumber\\
&\quad\quad\quad+
\frac{1}{2}\|\bar{x}^r-x^r\|^2+c\sqrt{K}\|x^{r+1}-x^r\|\left(\|x^{r+1}-x^r\|
+\|\bar{x}^{r}-x^r\|\right)\nonumber\\
&\le \zeta\|x^{r+1}-x^r\|^2+\zeta^{'}\|\bar{x}^r-x^r\|^2,\quad \mbox{for some }\zeta,\ \zeta'>0, \nonumber
\end{align}
where the first inequality follows from \eqref{eqBoundDifference}
and \eqref{eq:error_norm_AGP}, the second inequality is due to the Cauchy-Schwartz inequality and
the Lipschitz continuity of $\nabla \bar L_x(x;y^r)$, while the third inequality follows from the
fact that $\tilde{x}$ lies in the line segment joining $x^{r+1}$ and
$\bar{x}^r$ so that $\|\tilde{x}-x^r\|\le
\|x^{r+1}-x^r\|+\|\bar{x}^r-x^r\|$. This completes the proof. \QED

We then bound the decrease of the dual optimality gap.
\begin{lemma}\label{lm:dual-gap}
For each $r\ge 1$, there holds
\begin{equation}
\Delta_d^r-\Delta_d^{r-1}\le -\alpha(Ex^{r}-q)^T(E\barx^{r}-q). \label{eq:dual-gap}
\end{equation}
\end{lemma}
\proof
The reduction of the optimality gap in the dual space can be bounded as follows:
\begin{eqnarray*}
\Delta_d^r-\Delta_d^{r-1}&=&[d^*-d(y^{r})]-[d^*-d(y^{r-1})]\nonumber\\
&=& d(y^{r-1})-d(y^{r})\nonumber\\
&=&L(\barx^{r-1};y^{r-1})-L(\barx^{r};y^{r})\nonumber\\
&=&[L(\barx^{r};y^{r-1})-L(\barx^{r};y^{r})] + [L(\barx^{r-1};y^{r-1})-L(\barx^{r};y^{r-1})]\nonumber \\
&=&(y^{r-1}-y^{r})^T(q-E\barx^{r})+[L(\barx^{r-1};y^{r-1})-L(\barx^{r};y^{r-1})]\nonumber\\
&= &-\alpha(Ex^{r}-q)^T(E\barx^{r}-q)+[L(\barx^{r-1};y^{r-1})-L(\barx^{r};y^{r-1})]\nonumber\\
&\le &-\alpha(Ex^{r}-q)^T(E\barx^{r}-q), \quad \forall\; r\ge 1,
\end{eqnarray*}
where the last equality follows from the update of the dual variable
$y^{r-1}$, and the fact that $\bar{x}^{r-1}$ minimizes $L(\cdot,
y^{r-1})$. \QED

Lemma~\ref{lm:dual-gap} implies that if $q-Ex^r$ is close to the true dual gradient $\nabla d(y^r)=q-E\barx^r$, then the dual optimal gap is reduced after each ADMM iteration. However, since ADMM updates the primal variable by only one Gauss-Seidel sweep, the primal iterate $x^r$ is not necessarily close the minimizer $\barx^r$ of $L(x;y^r)$. Thus, unlike the method of multipliers (for which $x^r=\barx^r$ for all $r$), there is no guarantee that the dual optimality gap $\Delta_d^r$ is indeed reduced after each iteration of ADMM.

Next we proceed to bound the decrease in the primal gap $\Delta_p^r$.
\begin{lemma}\label{lm:primal-descent}
Suppose assumptions A(b) and A(f) hold. Then for each $r\ge1$, we
have
\begin{equation}
\Delta_p^{r}-\Delta_p^{r-1}\le \alpha\|Ex^r-q\|^2-\gamma\|x^{r+1}-x^r\|^2-\alpha(Ex^{r}-q)^T(E\barx^{r}-q)\label{eq:primal-gap}
\end{equation}
for some $\gamma$ independent of $y^r$.
\end{lemma}
\proof
Fix any $r\ge1$, we have
\[
L(x^r;y^{r-1})=f(x^r)+\langle y^{r-1}, q-Ex^r\rangle+\frac{\rho}{2}\|Ex^r-q\|^2\\
\]
and
\[
L(x^{r+1};y^r)=f(x^{r+1})+\langle y^r,q-Ex^{r+1}\rangle+\frac{\rho}{2}\|Ex^{r+1}-q\|^2.
\]
By the update rule of $y^r$ (cf.\ \eqref{eq:admm}), we have
\[
L(x^{r};y^r)=f(x^{r})+\langle y^{r-1}, q-Ex^{r}\rangle +\frac{\rho}{2}\|Ex^{r}-q\|^2+\alpha\|Ex^r-q\|^2.
\]
This implies
\[
L(x^r;y^r)=L(x^r;y^{r-1})+\alpha\|Ex^r-q\|^2.
\]
Recall from Lemma~\ref{lm:p-descent} that the alternating minimization of the Lagrangian function gives a sufficient descent. In particular, we have
\[
L(x^{r+1};y^r)-L(x^r;y^r)\le -\gamma\|x^{r+1}-x^r\|^2,
\]
for some $\gamma>0$ that is independent of $r$ and $y^r$. Therefore, we have
\[
L(x^{r+1};y^r)-L(x^r;y^{r-1})\le \alpha\|Ex^r-q\|^2-\gamma\|x^{r+1}-x^r\|^2, \quad \forall\; r\ge1.
\]
Hence, we have the following bound on the reduction of primal optimality gap
\begin{eqnarray*}
\Delta_p^{r}-\Delta_p^{r-1}&=&[L(x^{r+1};y^r)-d(y^r)]-[L(x^r;y^{r-1})-d(y^{r-1})]\nonumber\\
&=& [L(x^{r+1};y^r)-L(x^r;y^{r-1})]-[d(y^r)-d(y^{r-1})]\nonumber\\
&\le & \alpha\|Ex^r-q\|^2-\gamma\|x^{r+1}-x^r\|^2-\alpha(Ex^{r}-q)^T(E\barx^{r}-q), \quad \forall \; r\ge1,
\end{eqnarray*}
where the last step is due to Lemma~\ref{lm:dual-gap}. \QED

Notice that when $\alpha=0$ (i.e., no dual update in the ADMM algorithm), Lemma~\ref{lm:primal-descent} reduces to the sufficient decrease estimate \eqref{eq:p-descent} in Lemma~\ref{lm:p-descent}. When $\alpha>0$, the primal optimality gap is not necessarily reduced after each ADMM iteration due to the positive term $\alpha\|Ex^r-q\|^2$ in \eqref{eq:primal-gap}.
Thus, in general, we cannot guarantee a consistent decrease of either the dual optimality gap $\Delta_d^r$ or the primal optimality gap $\Delta_p^r$. However, somewhat surprisingly, the sum of the primal and dual optimality gaps decreases for all $r$, as long as the dual step size $\alpha$ is sufficiently small. This is used to establish the linear convergence of ADMM method.
\begin{theorem}\label{thm:main}
Suppose the conditions in Assumption A hold. 
Then 
the sequence of iterates $\{(x^r, y^r)\}$  generated by the ADMM
algorithm \eqref{eq:admm} converges linearly to an optimal
primal-dual solution for \eqref{eq:1}, provided the stepsize
$\alpha$ is sufficiently small. Moreover, the sequence of
feasibility violation $\{\|Ex^r-q\|\}$ also converges linearly.
\end{theorem}
\proof 
We show by induction that the sum of optimality gaps $\Delta^r_d+\Delta^r_p$ is reduced after each ADMM iteration, as long as the stepsize $\alpha$ is chosen sufficiently small. For any $r\ge1$, we denote
\begin{equation}\label{eq:useful}
\barx^r=\argmin_{\barx\in X(y^r)}\|\barx-x^r\|.
\end{equation}
By induction, suppose $\Delta^{r-1}_d+\Delta^{r-1}_p\le
\Delta^0_d+\Delta^0_p$ for some $r\ge 1$.  Recall that each $X_k$ is compact
and that the indicator function $i_{X_k}(x_k)$ is included in $h_k(x_k)$ (see the discussion after Assumption A), it follows that
$x^r\in X$, implying the boundedness of $x^r$. Thus, we obtain from Lemma
\ref{lm:eb} that
\begin{equation}\label{eq:bds}
\|x^r-\barx^r\|\le \tau \|\tilde \nabla L(x^r;y^r)\|
\end{equation}
for some $\tau>0$ (independent of $y^r$).
To prove Theorem~\ref{thm:main},
we combine the two estimates \eqref{eq:dual-gap} and \eqref{eq:primal-gap} to obtain
\begin{eqnarray}
[\Delta_p^{r}+\Delta_d^{r}]-[\Delta_p^{r-1}+\Delta_d^{r-1}]&=&
[\Delta_p^{r}-\Delta_p^{r-1}]+[\Delta_d^{r}-\Delta_d^{r-1}]\nonumber\\
&\le & \alpha\|Ex^r-q\|^2-\gamma\|x^{r+1}-x^r\|^2-2\alpha(Ex^{r}-q)^T(E\barx^{r}-q)\nonumber\\
&=&\alpha\|Ex^r-E\barx^r\|^2-\alpha\|E\barx^r-q\|^2-\gamma\|x^{r+1}-x^r\|^2.\label{eq:estimate}
\end{eqnarray}
Now we invoke \eqref{eq:bds} and Lemma~\ref{lm:estimate} to lower bound $\|x^{r+1}-x^r\|$:
\begin{equation}\label{eq:nice}
\|x^r-\barx^r\|\le\tau\|\tilde \nabla L(x^r;y^r)\|\le\tau \sigma \|x^{r+1}-x^r\|.
\end{equation}
Substituting this bound into \eqref{eq:estimate} yields
\begin{equation}\label{eq:descent}
[\Delta_p^{r}+\Delta_d^{r}]-[\Delta_p^{r-1}+\Delta_d^{r-1}] \le
(\alpha\|E\|^2\tau^{2}\sigma^{2}-\gamma)\|x^{r+1}-x^r\|^2-\alpha\|E\barx^r-q\|^2.
\end{equation}
Thus, if we choose the stepsize $\alpha$ sufficiently small so that
\begin{equation}\label{eq:alpha}
0<\alpha<\gamma\tau^{-2}\sigma^{-2}\|E\|^{-2},
\end{equation}
then the above estimate shows that
\begin{equation}\label{eq:induction}
[\Delta_p^{r}+\Delta_d^{r}]\le [\Delta_p^{r-1}+\Delta_d^{r-1}],
\end{equation}
which completes the induction. Moreover, the induction argument shows that if the stepsize $\alpha$ satisfies the condition \eqref{eq:alpha}, then the descent condition \eqref{eq:descent} holds for all $r\ge1$.

By the descent estimate \eqref{eq:descent}, we have
\begin{equation}\label{eq:limit}
\|x^{r+1}-x^r\|\to 0,\quad \|\nabla d(y^r)\|=\|E\barx^r-q\|\to 0.
\end{equation}


We now show that the sum of optimality gaps $\Delta^r_d+\Delta^r_p$
in fact contracts geometrically after a finite number of ADMM
iterations. By \eqref{eq:limit}, for any $\delta>0$, there must
exist a finite integer $\bar{r}>0$ such that for all $r\ge
\bar{r}$, $\|\nabla d(y^r)\|\le\delta$. Since $\Delta_d^r$, $\Delta^r_p$ are nonnegative and bounded from
above (see \eqref{eq:induction}), it follows that $d(y^r)$ is bounded from below by a constant $\zeta$ independent of $r$.
Applying the second part of Lemma~\ref{lm:eb}, we have that for all $r\ge \bar{r}$, the dual error
bound ${\rm dist}(y^r, Y^*)\le \tau \|\nabla d(y^r)\|$ holds true.

Therefore, it follows from
Lemma~\ref{lm:cost-to-go} that we have the following cost-to-go
estimate
\begin{equation}\label{eq:11}
\Delta_d^r=d^*-d(y^r)
\le  {\tau'}\|\nabla d(y^r)\|^2=  {\tau'}\|E\barx^r-q\|^2,
\end{equation}
for some $\tau'>0$ and for all $r\ge\bar{r}$.

Moreover, we can use Lemma~\ref{lm:cost-to-go} to bound
$\|x^{r+1}-x^{r}\|^2$ from below by $\Delta^r_p$. In particular, we
have from \eqref{eq:nice} and Lemma~\ref{lm:cost-to-go}  that
\begin{align*}
\Delta^r_p&\le \zeta\|x^{r+1}-x^r\|^2+\zeta^{'}\|\bar{x}^r-x^r\|^2\\
&\le\zeta\|x^{r+1}-x^r\|^2+\zeta^{'}\tau^2\sigma^2\|{x}^{r+1}-x^r\|^2\\
&=\left(\zeta+\zeta^{'}\tau^2\sigma^2\right)\|{x}^{r+1}-x^r\|^2.
\end{align*}
Substituting this bound and \eqref{eq:11} into \eqref{eq:descent},
and assuming that $\alpha>0$ satisfies \eqref{eq:alpha}, we obtain
\begin{eqnarray*}
[\Delta_p^{r}+\Delta_d^{r}]-[\Delta_p^{r-1}+\Delta_d^{r-1}]&\le &
(\alpha\|E\|^2\tau^2\sigma^2-\gamma)\|x^{r+1}-x^r\|^2-\alpha\|E\barx^r-q\|^2\\
&\le & -\frac{(\gamma-\alpha\|E\|^2\tau^2\sigma^2)}{\zeta+\zeta^{'}\tau^2\sigma^2}\Delta_p^r-\alpha(\tau')^{-1}\Delta_d^r\\
&\le & -{
\min\left\{\frac{(\gamma-\alpha\|E\|^2\tau^2\sigma^2)}{\zeta+\zeta^{'}\tau^2\sigma^2},\alpha(\tau')^{-1}\right\}}[\Delta_p^r
+\Delta_d^r].
\end{eqnarray*}
Since $\alpha>0$ is chosen small enough such that \eqref{eq:alpha}
holds, we have
\[
\lambda:=\min\left\{\frac{\gamma-\alpha\|E\|^2\tau^2\sigma^2}{\zeta+\zeta^{'}\tau^2\sigma^2},\alpha(\tau')^{-1}\right\}>0.
\]
Consequently, we have
\[
[\Delta_p^{r}+\Delta_d^{r}]-[\Delta_p^{r-1}+\Delta_d^{r-1}]\le -\lambda[\Delta_p^{r}+\Delta_d^{r}]
\]
which further implies
\[
0\le [\Delta_p^{r}+\Delta_d^{r}]\le \frac{1}{1+\lambda}[\Delta_p^{r-1}+\Delta_d^{r-1}].
\]
This shows that the sequence $\{\Delta_p^{r}+\Delta_d^{r}\}_{r\ge
\bar{r}}$ converges to zero Q-linearly\footnote{A sequence
$\{x^r\}$ is said to converge $Q$-linearly to some  $\bar{x}$ if
$\|x^{r+1}-\bar{x}\|/\|x^{r}-\bar{x}\|\le \mu$ for all $r$, where
$\mu\in(0,1)$ is some constant. A sequence $\{x^r\}$ is said to converge to $\bar{x}$
$R$-linearly if $\|x^r-\bar{x}\|\le c\mu^r$ for all $r$ and for some $c>0$. }.
As a result, we conclude that
$\{\Delta_p^{r}+\Delta_d^{r}\}$ and hence both $\Delta_p^r$ and
$\Delta^r_d$ globally converge to zero R-linearly\footnote{To see
that such R-linear convergence is in fact global, note that
$\bar{r}>0$ is finite, and $\Delta_p^{r}+\Delta_d^{r}$ is Q-linearly
convergent for $r\ge \bar{r}$. Then one can always find an appropriate constant $c$
such that
$\Delta_p^{r}+\Delta_d^{r}\le c(1+\lambda)^{-r}$ for all $r=1, 2,\ldots.$}.

We next show that the dual sequence $\{y^r\}$ 
is also R-linearly convergent. To this end,
notice that the inequalities \eqref{eq:nice} and \eqref{eq:descent}
imply
\begin{equation}\label{eq:descent1}
[\Delta_p^{r}+\Delta_d^{r}]-[\Delta_p^{r-1}+\Delta_d^{r-1}] \le
(\alpha\|E\|^2-\gamma\tau^{-2}\sigma^{-2})\|x^r-\barx^r\|^2-\alpha\|E\barx^r-q\|^2.
\end{equation}
Then by \eqref{eq:descent1}, we see that both $\|x^r-\barx^r\|\to 0$
and $\|E\barx^r-q\|\to 0$ R-linearly. This implies that $Ex^r-q\to
0$ R-linearly and $\nabla d(y^r)\to 0$ R-linearly. Using the fact
that ${\rm dist}(y^r,Y^*)\le \tau\|\nabla d(y^r)\|$, we conclude
that $y^r$ converges R-linearly to an optimal dual solution.

We now argue that the primal iterates $\{x^r\}$ converge to an optimal solution of \eqref{eq:1}.
By the inequality
\eqref{eq:descent}, we can further conclude that
\[
\|x^{r+1}-x^r\|^2\to0,\quad \|E\barx^r-q\|\to 0
\]
R-linearly. Notice that the R-linear convergence of
$\|x^{r+1}-x^r\|^2\to0$ implies that $\|x^{r+1}-x^r\|\to0$
R-linearly. This further shows that $x^r\to x^\infty$ R-linearly for
some $x^\infty$. Denote the limit of dual sequence $\{y^r\}$ by
$y^\infty$. By the preceding argument, we know $y^\infty$ is a dual
optimal solution of \eqref{eq:1}. To show that $x^\infty$ is a
primal optimal solution of \eqref{eq:1}, it suffices to prove that
$x^\infty\in X(y^\infty)$. Using \eqref{eq:nice}, and the fact that
$\|x^r-\bar{x}^r\|\to 0$, we have
$$\|{x}^{\infty}-\bar{x}^r\|\le\|x^r-{x}^{\infty}\|+\|x^r-\bar{x}^r\|\to 0 .$$
Since $\bar{x}^r\in X(y^r)$, we have $L(\bar{x}^r,y^r)\le L(x, y^r)$
for all $x\in X$. Passing limit, we obtain
$L({x}^{\infty},{y}^{\infty})\le L(x, {y}^{\infty})$ for all $x\in
X$, that is, ${x}^{\infty}\in X({y}^{\infty})$. It then follows that
the sequence $\{x^r\}$ converges R-linearly to a primal optimal
solution.
%
%
%
%
%
%
%
\QED

The following corollary relaxes the compactness assumption of the feasible set $X$ (Assumption A-(g) in Theorem~\ref{thm:main}),
and replaces it by the boundedness of the primal-dual iterates.

\begin{corollary}\label{co:1}
Suppose assumptions A(a)---A(f) hold, and that $X$ is a polyhedral
set. Further assume that either one of the following two assumptions
holds true
\begin{enumerate}
\item The sequence of dual iterates $\{y^r\}$ lies in a compact set,
and that the set $\{x \mid L(x,y)\le \zeta\}$ is compact for any
finite $y$ and $\zeta$.
\item The sequence of primal-dual iterates $\{(x^r, y^r)\}$ lies in a
compact set.
\end{enumerate}
Then if $\alpha$ is chosen sufficiently small $($cf.
\eqref{eq:alpha}$)$, all the conclusions stated in Theorem
\ref{thm:main} still hold true. Moreover, the sequence of function
values $\{f(x^r)\}$ also converges linearly.
\end{corollary}

\proof Suppose the first assumption is true, so that all the
dual iterates $\{y^r\}$ lie in a compact set $\tilde{Y}$. Let us define
$$\delta:=\sup\left \{\|x\|\mid [L(x,y)-d(y)]+[d^*-d(y)]\le \Delta^0_p+\Delta^0_p, \forall~y\in\tilde{Y}\right\}.$$

Clearly we have $[d^*-d(y)]\ge 0$ and $[L(x,y)-d(y)]\ge 0$ for all
$y\in\tilde{Y}$. Moreover $\delta<\infty$ due to the compactness of
$\tilde{Y}$ as well as the compactness of the set $\{x \mid
L(x,y)\le \zeta\}$ for any finite $y$.

Again we show by induction that the sum of optimality gaps
$\Delta^r_d+\Delta^r_p$ is reduced after each ADMM iteration, as
long as the stepsize $\alpha$ is chosen sufficiently small. For any
$r\ge 1$, by induction we suppose $\Delta^{r-1}_d+\Delta^{r-1}_p\le
\Delta^0_d+\Delta^0_p$. Then $\|x^r\|\le \delta<\infty$, and by the
first part of Lemma \ref{lm:eb}, the primal error bound holds. Then
we can carry out exactly the same analysis as the proof of Theorem
\ref{thm:main} to arrive at the same conclusion.

Additionally, since
\begin{eqnarray*}
f(x^{r+1})-d^*&=&[f(x^{r+1})-d(y^r)]+[d(y^r)-d^*]\\
&=&[L(x^{r+1};y^r)-d(y^r)]-[d^*-d(y^r)]-\langle y^r,q-Ex^{r+1}\rangle -\frac{\rho}{2}\|Ex^{r+1}-q\|^2\\
&=&\Delta_p^r-\Delta_d^r-\langle y^r,q-Ex^{r+1}\rangle -\frac{\rho}{2}\|Ex^{r+1}-q\|^2
\end{eqnarray*}
and
\[
\Delta_p^r\to 0,\quad \Delta_d^r\to 0,\quad \|q-Ex^r\|\to0,
\]
linearly, it follows that $f(x^{r+1})-d^*\to 0$ R-linearly (recall
that $y^r$ is now in a compact set). The proof is complete.

Similarly, when the second assumption is true, then  the error bound condition
in Lemma \ref{lm:eb} again holds, and by using the same argument
above we arrive at the desired conclusion.
\QED

As a remark, we point out that the proof of Corollary~\ref{co:1} also shows that the same linear convergence of
$f(x^r)\to d^*$ also holds under the assumptions in Theorem~\ref{thm:main}.


We close this section by providing a few examples that satisfy the
assumptions in Theorem \ref{thm:main} and Corollary
\ref{co:1}. First consider the following $\ell_1$ minimization
problem
\begin{align}
\min_{x} \|x\|_1, \quad {\rm s.t.} \ Ex=b, \ a\le x_k\le b, \
k=1,\cdots, K
\end{align}
which can be equivalently written as a $K$-block problem
\begin{align}
\min_{\{x_k\}} \sum_{k=1}^{K} |x_k|, \quad {\rm s.t.} \
\sum_{k=1}^{K}e_k x_k=b, \ a\le x_k\le b, \ k=1,\cdots, K,
\end{align}
where $e_k$ is the $k$-th column of $E$, $a$ and $b$ are some
scalars. It is easy to verify that this problem meets all the
conditions listed in Assumption A, hence the linear convergence result in
Theorem~\ref{thm:main} applies. The same is true for the following mixed
$\ell_1/\ell_2$ minimization problem
\begin{align}
\min_{\{x_k\}} \sum_{k=1}^{K}\|x_k\|, \quad {\rm s.t.}\
\sum_{k=1}^{K}E_k x_k=b,\ a\le x_k\le b, \ k=1,\cdots, K,
\end{align}
where $x_k$, $a$ and $b$  are now $n$-dimensional vectors, and
$E_k\in\R^{m\times n}$ is some matrix with full column rank.

Furthermore, the boundedness assumption in Corollary \ref{co:1} is satisfied
by many examples of the two-block ADMM described in
\cite{boyd}. This is because when $K=2$,
$\alpha/\beta\in(0,\frac{1}{2}(1+\sqrt{5}))$ and $E_k$'s are full
column rank, it is known that both the primal and dual iterates
generated by the two-block ADMM algorithm indeed lie in a bounded
set \cite{Glow84}. Therefore the second assumption made in the
Corollary \ref{co:1} holds true, hence we only require assumptions
A(a)--A(f), which are in fact quite mild. For example, they are met
by the following instance of the consensus problem (see
\cite[Section 7]{boyd} for introduction of the consensus problem)
\begin{align}\label{consensus}
\min_{\{x_k\},z}\sum_{k=1,\cdots, K}\|Ax_k-b\|^2+w\|x_k\|_1\quad{\rm
s.t.} \ x_k-z=0, \ k=1,\cdots, K,
\end{align}
where $w>0$ is some constant. Thus, the two block ADMM algorithm converges
linearly for \eqref{consensus} regardless of the rank of $A$. Note that when $A$ is full column
rank, the objective is strongly convex. Consequently the error bound
condition in Lemma \ref{lm:eb} holds true globally, and the
coefficient $\tau$ can be at least greater than $\lambda_{\rm
min}(A^T A)$.


%
%


\newsection{Variants of ADMM}\label{secVariants}

The convergence analysis of Section 3 can be extended to some variants of the ADMM. We briefly describe two of them below.

\subsection{Proximal ADMM}\label{subProximal}

In the original ADMM \eqref{eq:admm}, each block $x_k$ is updated by solving a convex optimization subproblem \emph{exactly}. For large scale problems, this subproblem may not be easy to solve unless the matrix $E_k$ is unitary (i.e., $E^T_kE_k=I$) in which case the variables in $x_k$ can be further decoupled (assuming $f_k$ is separable). If the matrix $E_k$ is not unitary, we can still employ a simple proximal gradient step to \emph{inexactly} minimize $L(x^{r+1}_1,...,x^{r+1}_{k-1},x_k,x^r_{k+1},...,x^r_K)$. More specifically, we update each block of $x_k$ according to the following procedure
\begin{eqnarray}
    x_k^{r+1}\!\!& =&\!\! {\rm arg}\!\min_{x_k}\Big\{\;h_k(x_k)+\langle y^r,q-E_kx_k\rangle+\left\langle A_k^T\nabla g_k(A_kx^r_k),x_k-x^r_k\right\rangle +\frac{\beta}{2}\|x_k-x^r_k\|^2\nonumber\\[10pt]
    &&\!\!+\Big\langle \rho E_k^T\Big(\sum_{j<k}E_jx^{r+1}_j+\sum_{j\ge k}E_jx_j^r-q\Big),x_k-x^r_k\Big\rangle\Big\} \label{eq:linearized-admm}
    \end{eqnarray}
in which the smooth part of the objective function in the $k$-th subproblem, namely,
$$
g_k(A_kx_k)+\langle y^r,q-E_kx_k\rangle+
\frac{\rho}{2}\Big\|E_kx_k+\sum_{j< k}E_jx_j^{r+1}+\sum_{j>
k}E_jx_j^{r}-q\Big\|^2
$$
is linearized locally at $x^r_k$, and a proximal term
$\frac{\beta}{2}\|x_k-x_k^r\|^2$ is added. Here, $\beta>0$ is a
positive constant. With this change, updating $x_k$ is easy when
$h_k$ (the nonsmooth part of $f_k$) is separable. For example, this
is the case for compressive sensing applications where
$h_k(x_k)=\|x_k\|_1$, and the resulting subproblem admits a closed
form solution given by the component-wise soft thresholding (also
known as the shrinkage operator). We note that the proximal ADMM
algorithm described here is slightly more general than the proximal
ADMM algorithm seen in the literature, in which only the
penalization term $\frac{\rho}{2}\Big\|E_kx_k+\sum_{j<
k}E_jx_j^{r+1}+\sum_{j> k}E_jx_j^{r}-q\Big\|^2 $ is linearized
locally at $x^r_k$; see e.g., \cite{YangZhang2011, WangYuan2012}.

We claim that Theorem~\ref{thm:main} holds for the proximal ADMM
algorithm {\it without requiring assumption A-(f)} (the full rankness
of $E_k$'s). Indeed, to establish the (linear) convergence of the
proximal ADMM \eqref{eq:linearized-admm}, we can follow the same
proof steps as that for Theorem~\ref{thm:main}, with the only
changes being in the proof of
Lemmas~\ref{lm:p-descent}-\ref{lm:estimate} and Lemma
\ref{lm:cost-to-go}. We first show that Lemma~\ref{lm:p-descent}
holds without assumption A(f). Clearly subproblem
\eqref{eq:linearized-admm} is now {\it strongly convex} without the
full column rank assumption of $E_k$'s made in A(f). In the
following, we will show that as long as $\beta$ is large enough,
there is a sufficient descent:
\begin{equation}\label{eq:s-descent}
{\color{black}
L(\vx^{r+1};y^r)-L(\vx^{r};y^r)\le -\gamma\|\vx^{r+1} - \vx^{r}\|^2,\quad \mbox{for some $\gamma > 0$ independet of $y^r$.}
}
\end{equation}
This property can be seen by bounding the smooth part of $L(x^{r+1}_1,...,x^{r+1}_{k-1},x_k,x^r_{k+1},...,x^r_K)$, which is given by
\[
{\bar L}_k(x_k):=g_k(A_kx_k)+\langle y^r,q-E_kx_k\rangle +\frac{\rho}{2}\Big\|\sum_{j<k}E_jx^{r+1}_j+\sum_{j> k}E_jx_j^r+E_kx_k-q\Big\|^2,
\]
with the Taylor expansion at $x^r_k$:
\begin{equation}\label{eq:estimate1}
{\bar L}_k(x^{r+1}_k)\le {\bar L}_k(x^r_k)+\langle \nabla {\bar L}_k(x^r_k),x_k^{r+1}-x^r_k\rangle +\frac{\nu}{2}\|x^{r+1}_k-x^r_k\|^2
\end{equation}
where
$$\nu:=L\|A_k\|\|A^T_k\|+\rho\|E_k^TE_k\|$$
is the Lipschitz constant of ${\bar L}_k(\cdot)$ and $L$ is the Lipschitz constant of $\nabla g_k(\cdot)$. Making the above inequality more explicit yields
\begin{eqnarray}
&&\!\!\!\!\!\!\!\!\!\!\!\!\!\!\!\!L(x^{r+1}_1,...,x^{r+1}_{k-1},x^{r+1}_k,x^r_{k+1},...,x^r_K;y^r)
-L(x^{r+1}_1,...,x^{r+1}_{k-1},x^r_k,x^r_{k+1},...,x^r_K;y^r)\nonumber
\\ [10pt]
&\le& h_k(\vx_k^{r+1})-h_k(\vx_k^r)+\langle y^r,E_k(x^r_k-x^{r+1}_k)\rangle+\left\langle A_k^T\nabla g_k(A_kx^r_k),x^{r+1}_k-x^r_k\right\rangle \nonumber\\
    &&\!\!+\left\langle \rho E_k^T\left(\sum_{j<k}E_jx^{r+1}_j+\sum_{j\ge k}E_jx_j^r-q\right),x^{r+1}_k-x^r_k\right\rangle +\frac{\nu}{2} \|\vx_k^{r+1}-\vx_k^r\|^2\nonumber\\ [10pt]
&\le& -\frac{\beta}{2}\|\vx_k^{r+1}-\vx_k^r\|^2 +\frac{\nu}{2}
\|\vx_k^{r+1}-\vx_k^r\|^2\nonumber\\ [10pt]
& =& -\gamma\|\vx^{r+1}_k-\vx^r_k\|^2, \quad \forall\; k,\label{eq:last-estimate}
\end{eqnarray}
provided the regularization parameter $\beta$ satisfies
$$
\gamma:=\frac{1}{2}\left(\beta-\nu\right)>0. 
$$
In the above derivation of \eqref{eq:last-estimate}, the first step is due to  \eqref{eq:estimate1}, 
while the second inequality follows from the definition of $x^{r+1}_k$ (cf.\ \eqref{eq:linearized-admm}). Summing \eqref{eq:last-estimate} over all $k$ yields the desired estimate of sufficient descent \eqref{eq:s-descent}.

To verify that Lemma~\ref{lm:estimate} still holds for the proximal
ADMM algorithm, we note from the corresponding optimality condition
for \eqref{eq:linearized-admm}
\[
x_k^{r+1}=\prox_{h_k}\left[x_k^{r+1}-A_k^T\nabla_{x_k}g_k(A_kx^{r}_k)+E_k^Ty^r-\rho E_k^T\left(\sum_{j< k}E_jx_j^{r+1}+\sum_{j\ge k}E_jx^r_j-q\right) -\beta(x_k^{r+1}-x_k^r)\right].
\]
Using this relation in place of \eqref{eq:old} and following the same proof steps, we can easily prove that the bound \eqref{eq:old-estimate} in Lemma~\ref{lm:estimate} can be extended to the proximal ADMM algorithm. Thus, the convergence results in Theorem~\ref{thm:main} remain true for the proximal ADMM algorithm \eqref{eq:linearized-admm}.

It remains to verify that Lemma~\ref{lm:cost-to-go} still holds
true. In fact the first part of Lemma~\ref{lm:cost-to-go} can be
shown to be independent of the iterates, thus it trivially holds
true for the proximal ADMM algorithm. To show that the second part
of Lemma~\ref{lm:cost-to-go} is true, note that the optimality
condition of the proximal ADMM algorithm implies that
\begin{align*}
x^{r+1}_k&={\rm
prox}_{h_k}\left[x^{r+1}_k-\nabla_{x_k}\bar{L}\left(\{x^{r+1}_{j<
k}\}, \{x^r_{j}\}_{j\ge k}; y^r\right)-\beta(x_k^{r+1}-x_k^r)\right]\\
&:= {\rm
prox}_{h_k}\left[x^r_k-\nabla_{x_k}\bar{L}(x^r;y^r)-e^r_k\right]
\end{align*}
where in this case $e^r_k$ is given as
\[
e^r_k:=x^r_k-x^{r+1}_k+\nabla_{x_k}\bar{L}\left(\{x^{r+1}_{j< k}\},
\{x^r_{j}\}_{j\ge k};
y^r\right)-\nabla_{x_k}\bar{L}(x^r;y^r)+\beta(x_k^{r+1}-x_k^r).\] It
is then straightforward to show that the norm of $e^r_k$ can be
bounded by $c^{'}\|x^r-x^{r+1}\|$ for some constant $c^{'}>0$. The
rest of the proof follows the same steps as in Lemma
\ref{lm:cost-to-go}.

\subsection{Jacobi Update}

Another popular variant of the ADMM algorithm is to use a Jacobi
iteration (instead of a Gauss-Seidel iteration) to update the primal
variable blocks $\{x_k\}$. In particular, the ADMM iteration
\eqref{eq:admm} is modified as follows:
\begin{equation}\label{eq:jacobi}
x^{r+1}_k=\argmin_{x_k}\left(h_k(x_k)+g_k(A_kx_k)-\langle y^r,E_kx_k\rangle+\frac{\rho}{2}\left\|E_kx_k+\sum_{j\neq k}E_jx_j^{r}-q\right\|^2\right),\quad \forall\; k.
\end{equation}

The convergence for this direct Jacobi scheme is unclear, 
as the augmented Lagrangian function may not decrease after each Jacobi update.
In the following, we consider a modified Jacobi
scheme with an explicit stepsize control. Specifically, let us introduce an intermediate variable $w=(w^T_1, \cdots,
w^T_K)^T\in \Re^n$. The modified Jacobi update is given as follows:
\begin{align}
w^{r+1}_k&=\argmin_{x_k}\left(h_k(x_k)+g_k(A_kx_k)-\langle
y^r,E_kx_k\rangle+\frac{\rho}{2}\left\|E_kx_k+\sum_{j\neq
k}E_jx_j^{r}-q\right\|^2\right),\quad \forall\;
k,\label{eq:jacobi_modif_w}\\
x^{r+1}_k&=x^r_k+\frac{1}{K}\left(w^{r+1}_k-x^r_k\right), \quad
\forall\; k. \label{eq:jacobi_modif_x}
\end{align}
where a stepsize of $1/K$ is used in the update of each variable block.

With this modification, we claim that
Lemmas~\ref{lm:p-descent}-\ref{lm:estimate} and Lemma
\ref{lm:cost-to-go} still hold.
In particular, Lemma~\ref{lm:p-descent} can be argued as follows.
The strong convexity of $L(x;y)$ with respect to the variable block
$x_k$ implies that
\begin{align}
&L\left(x^r_1, \cdots, x^r_{k-1}, x^r_k, x^r_{k+1}\cdots,x^r_K; y^r
\right)-L\left(x^r_1, \cdots, x^r_{k-1}, w^r_k,
x^r_{k+1}\cdots,x^r_K; y^r\right) \nonumber\\
&\ge \gamma\|w^{r+1}_k-x^r_k\|^2, \ \forall~k.\nonumber
\end{align}

Using this inequality we obtain
\begin{align*}
&L(x^r; y^r)-L(x^{r+1}; y^r)\\
&=L(x^r; y^r)-L\left(\frac{K-1}{K}x^r+\frac{1}{K}w^{r+1}; y^r\right)\nonumber\\
&= L(x^r; y^r)-L\left(\frac{1}{K}\sum_{k=1}^{K}(x^r_1, \cdots,
x^r_{k-1},
w^{r+1}_k, x^r_{k+1}\cdots,x^r_K); y^r\right)\\
&\ge L(x^r; y^r)-\frac{1}{K}\sum_{k=1}^{K}L\left(x^r_1, \cdots,
x^r_{k-1}, w^{r+1}_k, x^r_{k+1}\cdots,x^r_K; y^r\right)\\
&=\frac{1}{K}\sum_{k=1}^{K}\left(L(x^r; y^r)-L\left(x^r_1, \cdots,
x^r_{k-1}, w^{r+1}_k, x^r_{k+1}\cdots,x^r_K; y^r\right)\right)\\
&\ge\frac{\gamma}{K}\sum_{k=1}^{K}\|w^{r+1}_k-x^r_k\|^2\\
&=\frac{\gamma}{K}\|w^{r+1}-x^r\|^2.
\end{align*}
where the first inequality comes from the convexity of the augmented
Lagrangian function.

From the update rule \eqref{eq:jacobi_modif_x} we have
$K(x^{r+1}_k-x^r_k)=(w^{r+1}_k-x^r_k),$ which combined with the
previous inequality yields
\[L(x^r; y^r)-L(x^{r+1}; y^r)\ge \gamma K\|x^{r+1}-x^r\|^2.\]

The proof of Lemma~\ref{lm:estimate} also requires only minor
modifications. In particular, we have the following optimality
condition for \eqref{eq:jacobi}
\[
w_k^{r+1}=\prox_{h_k}\left[w_k^{r+1}-A_k^T\nabla_{x_k}g_k(A_kw^{r+1}_k)+E_k^Ty^r-\rho
E_k^T\left(\sum_{j\neq k}E_jx_j^{r}+E_kw^{r+1}_k-q\right) \right]
\]
Similar to the proof of Lemma~\ref{lm:estimate}, we have
\[
\left\|w_k^{r+1}-\prox_{h_k}\left[x^r_k-A_k^T\nabla_{x_k}g_k(A_kx^{r}_k)+E_k^Ty^r-\rho
E_k^T\left(E x^r-q\right) \right]\right\|\le c\|w^{r+1}-x^r\|.
\]
Utilizing the relationship $K(x^{r+1}_k-x^r_k)=(w^{r+1}_k-x^r_k)$, we can establish
Lemma~\ref{lm:estimate} by following similar proof steps (which we omit due to space reason).

Lemma \ref{lm:cost-to-go} can be shown as follows. We first
express $w_k^{r+1}$ as
\begin{align*}
w_k^{r+1}&=\prox_{h_k}\left[w_k^{r+1}-\nabla_{x_k}\bar{L}\left(\{x^{r}_{j\ne
k}\}, w^{r+1}_{k}; y^r\right)\right]\\
&=\prox_{h_k}\left[x_k^{r}-\nabla_{x_k}\bar{L}\left(x^{r};
y^r\right)-e^r_k\right]\\
\end{align*}
where we have defined
\[e^r_k:=\nabla_{x_k}\bar{L}\left(\{x^{r}_{j\ne
k}\}, w^{r+1}_{k}; y^r\right)-\nabla_{x_k}\bar{L}\left(x^{r};
y^r\right)+x^r_k-w_k^{r+1}.\] Again by using the relationship
$K(x^{r+1}_k-x^r_k)=(w^{r+1}_k-x^r_k)$, we can bound the norm of
$e^r_k$ by $c^{'}\|x^{r+1}-x^r\|$, for some $c^{'}>0$. The
remaining proof steps are similar to those in Lemma \ref{lm:cost-to-go}.

Since Lemmas~\ref{lm:p-descent}-\ref{lm:estimate} and Lemma
\ref{lm:cost-to-go} hold for the Jacobi version of the ADMM algorithm
with a step size control, we conclude that the convergence results of
Theorem~\ref{thm:main} remain true in this case.

\newsection{Concluding Remarks}

In this paper we have established the convergence and the rate of
convergence of the classical ADMM algorithm when the number of
variable blocks are more than two and in the absence of strong
convexity. Our analysis is a departure of the conventional analysis
of ADMM algorithm which relies on a contraction argument involving a weighted
(semi-)norm of $( x^r-x^*,y^r-y^*)$, see
\cite{Ga83,Gab76,Glow84,Glt89,HeTY12,HeYuan12,KoMe,TaoYuan}. In our
analysis, we require neither the strong convexity of the objective
function nor the row independence assumption of the constrained
matrix $E$. Instead, we use a local error bound to show that when
the stepsize of dual update is made sufficiently small, the sum of
the primal and the dual optimality gaps decreases after each ADMM
iteration, although separately they may individually increase. An
interesting issue for further research is to identify good practical stepsize
rules for dual update. While \eqref{eq:alpha} does suggest a dual  stepsize rule in terms of
error bound constants, it may be too conservative
and cumbersome to compute unless the objective function is strongly convex.
One possibility may be to use an adaptive dual stepsize rule to
guarantee the decrease of the sum of the primal and dual optimality
gaps.

\bigskip
\bigskip
\noindent {\bf Acknowledgement:} The authors are grateful to
Xiangfeng Wang and Dr.\ Min Tao of Nanjing University for their
constructive comments.


\newsection{Appendix}\label{secAppendix}

\subsection{Proof of Dual Error Bound \eqref{eq:dualeb}}\label{subProofDualEB}
The augmented Lagrangian dual function can be expressed as
\begin{align}
d(y)=\min_{x\in X}\langle
y,q-Ex\rangle+\frac{\rho}{2}\|q-Ex\|^2+g(Ax)+h(x).
\end{align}
For convenience, define $p(Ex):=\frac{\rho}{2}\|q-Ex\|^2$, and let
$\ell(x):=p(Ex)+g(Ax)+h(x)$. For simplicity, in this proof we
further restrict ourselves to the case where the nonsmooth part has
polyhedral level sets, i.e., $\{x: h(x)\le \xi\}$ is polyhedral for each $\xi$.
More general cases can be shown along similar lines, but the
arguments become more involved.

Let us define
$$x(y)\in \arg\min_{x\in X} \ell(x)+\langle y, q-Ex\rangle.$$

Let $(x^*, y^*)$ denote a primal and dual optimal
solution pair. Let $X^*$ and $Y^*$ denote the primal and dual optimal
solution set. 
The he following properties will be useful in our subsequent
analysis.
\begin{itemize}
\item[(a)] There exist positive scalars $\sigma_g$, $L_g$ such that $\forall~x(y),x(y')\in X$
\begin{enumerate}
\item[a-1)]$\langle A^T \nabla g(A x(y'))-A^T \nabla g(A x(y)),
x(y')-x(y)\rangle \ge \sigma_g\|A x(y')-A x(y)\|^2$.
\item[a-2)] $g(A x(y'))-g(Ax(y))-\langle A^T \nabla g(A x(y)),
x(y')-x(y)\rangle\ge \frac{\sigma_g}{2}\|A x(y')-Ax(y)\|^2.$
\item[a-3)] $\|A^T\nabla g(Ax(y'))-A^T\nabla g(A x(y))\|\le
L_g\|Ax(y')-Ax(y)\|.$
\end{enumerate}
\item[b)] All a-1)--a-3) are true for $p(\cdot)$ as well, with some
constants $\sigma_p$ and $L_p$.
\item[c)] $\nabla d(y)=q-Ex(y)$, and $\|\nabla d(y')-\nabla d(y)\|\le\frac{1}{\rho}\|y'-y\|.$
\end{itemize}

Part (a) is true due to the assumed Lipchitz continuity and strong
convexity of the function $g(\cdot)$. Part (b) is from the Lipchitz
continuity and strong convexity of the quadratic penalization
$p(\cdot)$. Part (c) has been shown in Lemma
\ref{lm:const-derivative} and Lemma \ref{lm:3.2}.

To proceed, let us rewrite the primal problem equivalently as
\begin{align}
d(y)=\min_{(x,s): x\in X, h(x)\le s} \langle y, q-Ex\rangle+
p(Ex)+g(Ax)+s.\label{eq:problemPolyhedral}
\end{align}
Let us write the polyhedral set $\{(x,s): x\in X, h(x)\le s\}$
compactly as $C_x x+ C_s s \ge c$ for some matrices
$C_x\in\mathbb{R}^{j\times n}$, $C_s\in\mathbb{R}^{j\times 1}$ and
$c\in\mathbb{R}^{j\times 1}$, where $j$ is some integer.
For any fixed $y$, let $(x(y), s(y))$ denote
one optimal solution for \eqref{eq:problemPolyhedral}, note we must
have $h(x(y))=s(y)$. Due to equivalence, if $y^*\in Y^*$, we must
also have $x(y^*)\in X^*$.

Define a set-valued function $\cM$ that assigns the vector $(d,
e)\in\mathbb{R}^{n}\times\mathbb{R}^{m}$ to the set of vectors
$(x,s,y,\lambda)\in\mathbb{R}^n\times\mathbb{R}\times
\mathbb{R}^m\times \mathbb{R}^j$ that satisfy the following system of
equations
\begin{align*}
E^Ty + 
 C_x^T\lambda &=d,\nonumber \\
C^T_s\lambda&=1,\\
q-Ex&=e,\\
\lambda\ge 0, \ (C_x x+C_s s)\ge c, \ \langle C_x x+C_s
s-c,\lambda\rangle&=0.
\end{align*}
It is easy to verify by using the optimality condition for problem
\eqref{eq:problemPolyhedral} that
\begin{align}
(x,s,y,\lambda)\in\cM(E^T\nabla p(Ex)+A^T\nabla g(Ax),
e)~\mbox{for some} \ \lambda\nonumber\\
\quad\quad\quad \mbox{if and only if}\  x=x(y), \ e=\nabla
d(y).\label{eq:PolyHedronOPT}
\end{align}
We can take $e=0$, and use the fact that $x(y^*)\in X^*$, we see
that $(x,s,y,\lambda)\in\cM(E^T\nabla p(Ex)+A^T\nabla g(Ax), 0)$ if
and only if $x\in X^*$ and $y\in Y^*$.

The following result states a well-known local upper Lipschitzian
continuity property for the polyhedral multifunction $\cM$; see
\cite{Hoffman, [16], luo-tseng}.

\begin{proposition}\label{propLip}
There exists a positive scalar $\theta$ that depends on $A, E, C_x,
C_s$ only, such that for each $(\bar{d}, \bar{e})$ there is a
positive scalar $\delta'$ satisfying
\begin{align}
\cM(d, e)\subseteq \cM (\bar{d}, \bar{e})+\theta
\|(d, e)-(\bar{d}, \bar{e})\|\mathcal{B},\\
 \quad\quad\quad \mbox{whenever} \
\|(d, e)-(\bar{d}, \bar{e})\|\le \delta'.
\end{align}
where $\mathcal{B}$ denotes the unit Euclidean ball in
$\mathbb{R}^n\times\mathbb{R}^{m}\times \mathbb{R}\times
\mathbb{R}^j$.
\end{proposition}

The following is the main result for this appendix. Note that the
scalar $\tau$ in the claim is independent the choice of $y$, $x$,
$s$, and is independent on the coefficients of the linear term $s$.
\begin{claim}
Suppose all the assumptions in Assumption A are satisfied. Then
there exits positive scalars $\delta$, $\tau$ such that for all
$y\in \cal U$ and $\|\nabla d(y)\|\le \delta$, there holds ${\rm
dist}(y, Y^*)\le \tau \|\nabla d(y)\|$.
\end{claim}
\proof By the previous claim, $\cM$ is locally Lipschitzian with
modulus $\theta$ at $(\nabla\ell(x^*), 0)=(E^T\nabla
p(Ex^*)+A^T\nabla g(Ax^*), 0)$.

Let $\delta\le \delta'/2$. We first show that if $\|\nabla d(y)\|\le
\delta$, then we must have $\|\nabla \ell(x(y))-\nabla
\ell(x^*)\|\le \delta'/2$. To this end, take a sequence
$y^1,y^2,\cdots,$ such that $e^r:=\nabla d(y^r)\to 0$. By assumption
A(g) $\{x(y^r)\}$ lies in a compact set. Due to the fact that
$s(y^r)=h(x(y^r))$, so the sequence $\{s(y^r)\}$ also lies in a compact set
(cf. Assumption A(e)). By passing to a subsequence if
necessary, let $(x^{\infty},s^{\infty})$ be a cluster
point of $\{x(y^r), s(y^r)\}$. In light of the continuity of $\nabla \ell
(\cdot)$, we have $(\nabla \ell(x(y^r)), e^r)\to
(\nabla\ell(x^{\infty}),0)$. Now for all $r$, $\{(x(y^r), s(y^r),
\nabla\ell(x(y^r)), e^r)\}$ lies in the set
$$\left\{(x, s, d, e)\mid (x,s,y,\lambda)\in\cM(d, e)~\mbox{for some}~ (y,\lambda)\right\}$$
which is polyhedral and thus is closed. Then we can pass limit to it
and conclude (cf. Proposition \ref{propLip})
$$(x^\infty, s^{\infty}, y^{\infty}, \lambda^{\infty})\in\cM (\nabla\ell(x^{\infty}), 0)$$
for some $(y^{\infty}, \lambda^{\infty})\in\mathbb{R}^{m}\times
\mathbb{R}^j$. Thus by \eqref{eq:PolyHedronOPT} and the discussions
that follow, we have $x^\infty\in X^*$ and $y^{\infty}\in Y^*$. By Lemma~\ref{lm:const-derivative}, we have
$\nabla \ell(x^*)=\nabla\ell(x^\infty)$, which further implies that $\nabla \ell(x(y^r))\to \nabla \ell(x^*)$. This
shows that the desired $\delta$ exists.

Then we let $e= \nabla d(y)$, and suppose $\|e\|\le \delta$. From
the previous argument we have
$$\|\nabla\ell(x(y))-\nabla\ell(x^*)\|+\|e\|\le \delta'/2+\delta'/2=\delta'.$$

Using the results in Proposition~\ref{propLip}, we have that there exists $(x^*,s^*, y^*, \lambda^*)\in \cM(\nabla \ell(x^*), 0)$ satisfying
\begin{align}
\|(x(y), s, y, \lambda)-(x^*,s^*, y^*, \lambda^*)\|&\le
\theta\left(\|\nabla\ell(x^*)-\nabla\ell(x(y))\|+\|e\|\right)\nonumber.
\end{align}

Since $(x(y),s,y,\lambda)\in\cM(\nabla \ell (x(y)), e)$, it follows from the
definition of $\cM$ that
\begin{align}
E^Ty 
+ C_x^T\lambda &=\nabla\ell(x(y)),\label{eq:FirstOrderPrimal1} \\
C^T_s\lambda&=1, \label{eq:FirstOrderPrimal2} \\
q-Ex(y)&=e,\label{eq:FeaibilityPrimal} \\
\lambda\ge 0, \ (C_x x(y)+C_s s(y))\ge c, \ \langle C_x x(y)+C_s
s(y)-c,\lambda\rangle&=0. \label{eq:ComplementarityPrimal}
\end{align}
Since $(x^*,s^*,y^*,\lambda^*)\in\cM(\nabla \ell (x^*), 0)$, we have from
the definition of $\cM$
\begin{align}
E^Ty^* + 
C_x^T\lambda^* &=\nabla\ell(x^*),\label{eq:FirstOrderOPT1} \\
C^T_s\lambda^*&=1, \label{eq:FirstOrderOPT2}\\
q-Ex^*&=0,\label{eq:FeaibilityOPT}\\
\lambda^*\ge 0, \ (C_x x^*+C_s s^*)\ge c, \ \langle C_x x^*+C_s
s^*-c,\lambda^*\rangle&=0.\label{eq:ComplementarityOPT}
\end{align}

Moreover, we have
\begin{align}
&\sigma_g\|A(x(y)-x^*)\|^2+\sigma_p\|E(x(y)-x^*)\|^2\nonumber\\
&\le\langle A^T\nabla g(Ax(y))-A^T\nabla g(Ax(y^*)),
x(y)-x(y^*)\rangle+\langle E^T\nabla p(Ex(y))-E^T\nabla p(Ex(y^*)),
x(y)-x(y^*)\rangle\nonumber\\
&=\langle\nabla\ell(x(y))-\nabla\ell(x(y^*)),
x(y)-x(y^*)\rangle\nonumber\\
&=\langle \lambda-\lambda^*, C_x x(y)-C_x x^*\rangle+\langle y-y^*,
E x(y)-E x^*\rangle
\nonumber
\end{align}
where the first inequality comes from the strong convexity of
$g(\cdot)$ and $p(\cdot)$; the last equality is from
\eqref{eq:FirstOrderPrimal1} and \eqref{eq:FirstOrderOPT1}.
Moreover, we have
\begin{align}
&\langle \lambda-\lambda^*, C_x x(y)-C_x
x^*\rangle\nonumber\\
&=\langle \lambda-\lambda^*, C_x x(y)-C_x
x^*\rangle+\langle\lambda-\lambda^*, C_s s-C_s s^*\rangle \nonumber\\
&=\langle \lambda-\lambda^*, (C_x x(y)+C_s s)-(C_x x^*+C_s
s^*)\rangle\nonumber\\
&=-\langle \lambda^*, C_x x(y)+C_s s-c\rangle-\langle \lambda, C_x
x^*+C_s s^*-c\rangle\le 0
\end{align}
where in the first equality we have used the fact that
$C^T_s\lambda-C^T_s\lambda^*=0$; see \eqref{eq:FirstOrderPrimal2}
\eqref{eq:FirstOrderOPT2}; in the third equality and in the last
inequality we have used the complementary conditions
\eqref{eq:ComplementarityOPT} and \eqref{eq:ComplementarityPrimal}.
As a result, we have
\begin{align}
&\sigma_g\|A(x(y)-x^*)\|^2+\sigma_p\|E(x(y)-x^*)\|^2\nonumber\\
&\le\langle y-y^*, (Ex(y)-q)-(Ex^*-q)\rangle\le \|y-y^*\|\|e\|,
\end{align}
where the last step is due to $\nabla d(y)=Ex(y)-q$ and $\nabla d(y^*)=Ex^*-q=0$.
Finally we have from Proposition~\ref{propLip}
\begin{align*}
&\|(x(y), s, y, \lambda)-(x^*,s^*, y^*, \lambda^*)\|^2\\
&\le
\theta^2\left(\|\nabla\ell(x^*)-\nabla\ell(x(y))\|+\|e\|\right)^2\nonumber\\
&\le\theta^2\left(2\|\nabla\ell(x^*)-\nabla\ell(x(y))\|^2+2\|e\|^2\right)\nonumber\\
&\le 2\theta^2\left(2\|\nabla g(x^*)-\nabla g(x(y))\|^2+2\|\nabla p(x^*)-\nabla p(x(y))\|^2+\|e\|^2\right)\nonumber\\
&\le2\theta^2\left(L^2_g\|A^T(x(y)-x^*)\|^2+L^2_p\|E^T(x(y)-x^*)\|^2+\|e\|^2\right)\nonumber\\
&\le
2\theta^2\max\left(\frac{2L^2_g}{\sigma_g},\frac{2L^2_p}{\sigma_p},1\right)\left(\sigma_g\|A^T(x(y)-x^*)\|^2+\sigma_p\|E^T(x(y)-x^*)\|^2+\|e\|^2\right)\\
&\le
2\theta^2\max\left(\frac{2L^2_g}{\sigma_g},\frac{2L^2_p}{\sigma_p},1\right)\left(\|e\|\|y-y^*\|+\|e\|^2\right)\nonumber\\
&\le
2\theta^2\max\left(\frac{2L^2_g}{\sigma_g},\frac{2L^2_p}{\sigma_p},1\right)\left(\|e\|\|(x(y),
s,y, \lambda)-(x^*, s^*, y^*, \lambda^*)\|
+\|e\|^2\right)\nonumber,
\end{align*}
where the second inequality is due to $\nabla \ell(x)=\nabla g(x)+\nabla p(x)$ and the fourth step follows from properties a-3) and b).

We see that the above inequality is quadratic in $\|(x(y), s, y,
\lambda)-(x^*, s^*, y^*, \lambda^*)\|/\|e\|$, so we have
$$\|(x(y),s, y, \lambda)-(x^*, s^*, y^*, \lambda^*)\|/\|e\|\le \tau$$ for some
scalar $\tau$ depending on $\theta$, $L_g$, $L_p$, $\sigma_g$,
$\sigma_p$.  It is worth noting that $\tau$ does not depend on the
choice of the coefficients of the linear term $s$. We conclude ${\rm
dist}(y, Y^*)\le \tau\|\nabla d(y)\|$.
 \QED

\vfill
\eject

\end{document}